\documentclass[11pt,a4paper]{article}
\usepackage[utf8]{inputenc}
\usepackage{lmodern}
\usepackage[T1]{fontenc}
\usepackage[english]{babel}
\usepackage{ifpdf}

\usepackage[left=1in, right=1in, top=1in, bottom=1in]{geometry}
\usepackage[dvipsnames]{xcolor}

\usepackage[square,numbers]{natbib}
\usepackage[colorlinks=true,linkcolor=teal,citecolor=teal]{hyperref}

\usepackage{graphicx, amsmath, amsthm, amssymb, enumitem, mathrsfs} 
\usepackage{subcaption}
\usepackage{booktabs}

\usepackage[capitalize]{cleveref}
\Crefname{subsection}{subsection}{subsections}
\Crefname{subsection}{Subsection}{Subsections}
\usepackage{soul}

\usepackage{tikz}
\usepackage{tikz-3dplot} 
\usepackage{pgfplots, pgfplotstable}
\pgfplotsset{compat=1.18}
\usetikzlibrary{3d} 
\usetikzlibrary{arrows.meta}

\usepackage{algorithm}
\usepackage{algorithmic}
\usepackage{authblk}

\newcommand{\Wass}{\mathcal{W}}
\newcommand{\Maas}{\mathcal{M}}
\newcommand{\Resis}{\mathcal{R}}

\newcommand{\Div}{\mathsf{Div}}

\newcommand{\Tree}{\mathcal{T}}
\newcommand{\Forest}{\mathcal{F}}

\newcommand{\Prob}{\mathscr{P}}

\newcommand{\Dist}{\Delta}

\newcommand{\Hitting}{\mathrm{H}}

\newcommand{\Deg}{\mathrm{d}}
\newcommand{\Vol}{\mathrm{Vol}}

\newcommand{\ep}[1]{\mathbb{E}\left({#1}\right)}
\newcommand{\pr}[1]{\mathbb{P}\left({#1}\right)}

\newtheorem{theorem}{Theorem}[section]
\newtheorem{corollary}[theorem]{Corollary}
\newtheorem{remark}[theorem]{Remark}
\newtheorem{definition}[theorem]{Definition}
\newtheorem{lemma}[theorem]{Lemma}
\newtheorem{proposition}[theorem]{Proposition}

\usepackage{xcolor,soul}

\title{Resistance Distance and \\ Linearized Optimal Transport on Graphs}

\hypersetup{
	pdftitle={Resistance Distance and Linearized OT on Graphs},
	pdfauthor={Robertson, Wan, Cloninger}
} 

\author{Sawyer Jack Robertson\footnote{Department of Mathematics, UC San Diego}\,, Zhengchao Wan\footnote{Department of Mathematics, University of Missouri}\,, and Alexander Cloninger${}^\ast{}$\footnote{Halicio\u{g}lu Data Science Institute, UC San Diego}}

\date{\today}

\begin{document}
	
\maketitle

\vspace*{-1cm}

\begin{abstract}
    We study the linearization of a discrete transportation distance between probability distributions on finite weighted graphs originally due to Maas (``Gradient flows of the entropy for finite {M}arkov chains,'' J. Funct. Anal. 261(8), 2011) which demonstrates various connections to the underlying combinatorial structure of the graph. For a connected graph and a reference density $\mu$ on its vertices, our main result is a nonasymptotic local linearization theorem showing that if $\nu$ is a small additive perturbation of $\mu$ then their squared discrete transportation distance is controlled above and below by the quadratic form of the pseudoinverse of a re-weighted graph Laplacian matrix. When the reference measure is stationary for the simple random walk on the graph, the weights agree with the original graph and this yields the quadratic form $(\mu-\nu)^\top L_w^\dagger (\mu-\nu)$, which can be viewed as a form of resistance distance between probability measures. This distance has a number of combinatorial and variational characterizations, including Beckmann and Benamou--Brenier formulas, a dual homogeneous Sobolev norm formula, a spanning $2$-forest formula, and a representation through random walk hitting times. Finally, we show that on the resulting ``resistance manifold,'' the gradient flow of the $\chi^2$ functional is the continuous-time random walk and that its geodesic strong convexity modulus equals the spectral gap of the normalized Laplacian. From this geometric vantage point, one recovers the classical fact that the spectral gap of the normalized Laplacian controls the exponential convergence rate of the random walk to stationarity.
\end{abstract}

\section{Introduction}
Recently many researchers have investigated \emph{linearization} as a means to circumvent computational bottlenecks associated with optimal transportation problems in large-scale data science applications. A number of methods have been developed along this line, generally according to the following rough template: on a metric space $\mathcal{X}$, equip a suitable class of probability measures $\Prob(\mathcal{X})$ with a transportation metric $\mathcal{T}$ so that $(\Prob(\mathcal{X}),\mathcal{T})$ has a (formal) Riemannian structure, then linearize $\mathcal{T}$ at a \emph{reference measure} $\mu\in\Prob(\mathcal{X})$ by identifying the Riemannian inner product on the tangent space at $\mu$ and using its induced norm distance as a surrogate. The quality of this linearization is usually governed by local geometry (e.g., curvature) of $(\Prob(\mathcal{X}),\mathcal{T})$ near $\mu$. 
 
As an example, consider the \emph{Otto calculus}~\cite{otto2001geometry} (the nonasymptotic bounds given below were proved in~\cite{peyre2018comparison}; see also~\cite[Sec. 5.5.2]{santambrogio2015optimal}). Let $\Omega\subseteq\mathbb{R}^{d}$ be a bounded convex domain, let $dx$ denote the normalized Lebesgue measure on $\Omega$, and let $h\in L^{1}(\Omega)$ satisfy $\int_{\Omega} h(x)\,dx = 0$. For $\varepsilon\in\mathbb{R}$ with $|\varepsilon|$ sufficiently small, assume the perturbed measure $dx+\varepsilon h(x)\,dx$ has density bounded from above and below, respectively, by constants $0 < a < b < \infty$, i.e.,
    \begin{align*}
        a \le 1+\varepsilon h(x) \le b \quad \text{for a.e. } x\in\Omega.
    \end{align*}
Then the $2$-Wasserstein distance, denoted $\Wass_{2}$, satisfies
    \begin{align}\label{eq:sobolev-linearization}
        b^{-1/2}\|\varepsilon h\|_{\dot{H}^{-1}(\Omega)} \leq \Wass_{2}(dx,\,dx+\varepsilon h(x)\,dx) \leq a^{-1/2} \|\varepsilon h\|_{\dot{H}^{-1}(\Omega)}.
    \end{align}
Here, $\|\cdot\|_{\dot{H}^{-1}(\Omega)}$ is the dual norm of the homogeneous Sobolev space $\dot{H}^{1}(\Omega)$ of mean-zero Lebesgue measurable functions with square-integrable weak gradients, equipped with the norm
    \begin{align*}
        \|f\|_{\dot{H}^1(\Omega)}^2 &= \int_{\Omega} |\nabla f(x)|^2\,dx,\quad f\in\dot{H}^{1}(\Omega).
    \end{align*}
Consequently,
    \begin{align}\label{eq:sobolev-linearization-asymptotic}
        \Wass_{2}(dx,\,dx+\varepsilon h(x)\,dx)^2 &= |\varepsilon|^2\,\|h\|_{\dot{H}^{-1}(\Omega)}^2 + o(|\varepsilon|^2).
    \end{align}
In the linearization framework above,~\cref{eq:sobolev-linearization} says that, after formally identifying tangent vectors at $dx$ with mean-zero perturbations $h$, the Riemannian metric tensor of $\Wass_2$ at $dx$ is the $\dot{H}^{-1}$ inner product. Hence, for small $|\varepsilon|$, $\Wass_2(dx,\,dx+\varepsilon h\,dx)$ is approximated by $|\varepsilon|\,\|h\|_{\dot{H}^{-1}(\Omega)}$.

We develop an analogous program on weighted graphs, establishing connections between a discrete transportation metric, its linearization, and various combinatorial properties of the underlying graph. We consider graphs of the form $G=(V,E,w)$, where $V=\{x_1,x_2,\dotsc, x_n\}$ is a finite set of vertices, $E\subseteq\binom{V}{2}$ is a collection of $m\geq 1$ undirected edges, and $w=(w_{xy})_{\{x, y\}\in E}$ is an assignment of positive edge weights. We denote by $\mathbb{R}^{V}$ the linear space of functions from $V$ to $\mathbb{R}$ identified with $\mathbb{R}^{n}$, and by $\Prob(V)$ the simplex of probability density functions on $V$. 

In our setting, the aforementioned linearization template is difficult to reconcile, mostly because of the metric choice. For distributions on $\mathbb{R}^{d}$, linearization approaches are often designed around the Monge formulation (see~\nameref{subsec:related-work}), but on discrete spaces, Monge transport is generally ill-posed. Its relaxation, the $p$-Wasserstein metric (see~\Cref{defn:wasserstein}), still exhibits geometric degeneracies: optimal couplings are typically non-unique, so the map $(\mu,\nu)\mapsto \Wass_p(\mu,\nu)$ is only piecewise smooth on the simplex, with kinks when the combinatorial type of an optimal plan changes. Geodesic behavior is also problematic. While it is true that $\Wass_1$ geodesics always exist for Polish metric spaces~\cite[Theorem 5.1]{bottou2018geometrical} (see also~\cite{hillion2014interpolation}), $\Wass_p$ geodesics do not exist in general for $p>1$, in particular, Maas~\cite{maas2011gradient} showed that even on a two-vertex graph, $\Wass_2$ may fail to admit unit-speed geodesics between given endpoint measures. 

Motivated by these geometric considerations, and inspired by the well known Benamou--Brenier formula for $\Wass_2$ (see~\cite{benamou2000computational}), Maas~\cite{maas2011gradient} introduced a ``Wasserstein-like'' transportation metric under which the interior of the probability simplex demonstrates a number of desirable geometric properties. We recall its definition after introducing some preliminary notation.

We denote by $\Prob_\star(V)$ the interior of $\Prob(V)$ (i.e., the set of everywhere positive densities). We denote by $A_{w}\in\mathbb{R}^{V\times V}$ the (weighted) \emph{adjacency matrix} of $G$, and for a vertex $x\in V$, we define its \emph{degree} $\Deg_{x}$ to be the sum of the weights of edges incident to $x$, i.e., $\Deg_{x} := \sum_{y\sim x} w_{xy}$. The \emph{degree matrix} $D_{w}\in\mathbb{R}^{V\times V}$ is the diagonal matrix of vertex degrees, and the \emph{Laplacian matrix} $L_{w}\in\mathbb{R}^{V\times V}$ is defined via the formula $L_{w} := D_{w} - A_{w}$. Recall that the \emph{simple random walk on $V$} is the Markov chain $(X_t)_{t\geq 0}$ with transition probability matrix $P=D_w^{-1} A_{w}$. We recall that if $G$ is connected then $(X_t)_{t\geq 0}$ is irreducible and has a unique stationary distribution $\pi\in\Prob(V)$, and if $G$ is additionally not bipartite, then $(X_t)_{t\geq 0}$ is also aperiodic and hence ergodic.

For scalars $a, b\geq 0$, we define their \emph{logarithmic mean} $\theta(a, b)$ via the formula $\theta(a, b) = \int_0^1 a^{1-\lambda} b^{\lambda} d\lambda$. Let $\mathcal{C}^1([0, 1];\, \Prob_\star(V))$ denote the set of $\mathcal{C}^{1}$ functions taking values in $\Prob_\star(V)$, and let $\mathcal{B}([0, 1];\, \mathbb{R}^{V})$ denote the set of Borel measurable functions from $[0, 1]$ to $\mathbb{R}^{V}$.

\begin{definition}[Discrete transportation distance~\cite{maas2011gradient}]\label{defn:mass-distance}
    Let $G=(V,E, w)$ be a connected graph and let $\pi\in\Prob_\star (V)$ denote the corresponding stationary distribution of the simple random walk on $V$. For any $(\psi_t)_{t\in[0,1]}\in \mathcal{B}([0, 1];\, \mathbb{R}^{V})$ and any $(\rho_t)_{t\in[0,1]}\in \mathcal{C}^{1}([0, 1];\, \Prob_\star (V))$, define the \emph{action functional}
        \begin{align}\label{eq:action-functional}
            \mathcal{A}((\psi_t), (\rho_t)) &:=  \int_0^{1} \sum_{\{x, y\}\in E} w_{xy}\, (\psi_t(x) - \psi_t(y))^2\, \theta\left(\frac{\rho_t(x)}{\pi(x)}, \frac{\rho_t(y)}{\pi(y)}\right)\, dt.
        \end{align}
    For $\mu,\nu\in\Prob_\star(V)$, their \emph{discrete transportation distance} $\Maas (\mu,\nu)$ is given by
        \begin{align}\label{eq:defn-maas-distance}
            \Maas (\mu,\nu)^2 &= \inf\Bigl\{ \mathcal{A}((\psi_t), (\rho_t)) \;:\; (\rho_t)_{t\in[0,1]}\in \mathcal{C}^{1}([0, 1];\,\Prob_\star(V)),\; (\psi_t)_{t\in[0,1]}\in \mathcal{B}([0, 1];\,\mathbb{R}^{V}),\notag\\
            &\qquad\qquad \frac{d}{dt}\rho_t(x) + \sum_{y\sim x} w_{xy}\, (\psi_t(y) - \psi_t(x))\, \theta\left(\frac{\rho_t(x)}{\pi(x)}, \frac{\rho_t(y)}{\pi(y)}\right) = 0,\notag\\
            &\qquad\qquad \rho_0 = \mu,\ \rho_1 = \nu \Bigr\}.
        \end{align}
\end{definition}

Regarding~\Cref{defn:mass-distance}, we remark that in~\cite{maas2011gradient} the author considered a broader class of means $\theta(\cdot,\cdot)$; we specialize to the logarithmic mean for simplicity and concreteness. Our conventions for density and action functional normalizations also differ slightly: the former amounts to a diagonal change of variables, and the latter is absorbed into constants as needed.

Maas~\cite[Theorem 1.1]{maas2011gradient} proves that $\Maas(\cdot,\cdot)$ is a metric on $\Prob_\star(V)$, and moreover that $(\Prob_\star(V),\Maas)$ defines a Riemannian manifold. Additionally, Maas shows that the gradient flow of the entropy functional on $\Prob_\star(V)$ with respect to the Riemannian metric tensor induced by $\Maas(\cdot,\cdot)$ coincides with the associated simple random walk diffusion. To our knowledge, however, no work has yet investigated the properties of linearized $\Maas(\cdot,\cdot)$ or the geometry of the resulting flat metric. 

\paragraph{Our main contribution.}

Our main result, which can be viewed as an analogue of~\cref{eq:sobolev-linearization}, asserts that near a fixed reference measure $\mu\in\Prob_\star(V)$, the local geometry of the discrete transportation distance is governed by a Laplacian associated with a $\mu$-dependent reweighting of the graph.

For a matrix $X$, we denote by $X^\dagger$ its Moore--Penrose inverse. We denote by $\mathbb{R}^{V}_0$ the subspace of $\mathbb{R}^{V}$ consisting of mean zero functions.

\begin{theorem}[Linearization of discrete transportation distance near a fixed measure]\label{thm:linearization}
    Let $G=(V,E, w)$ be a connected graph and let $\pi\in\Prob_\star(V)$ denote the stationary distribution of the simple random walk on $V$. Let $\mu\in\Prob_\star(V)$ be a fixed reference measure and define $w^{(\mu)}=(w^{(\mu)}_{xy})_{\{x,y\}\in E}$ according to the formula
        \begin{align*}
            w^{(\mu)}_{xy} := w_{xy}\,\theta\left(\frac{\mu(x)}{\pi(x)},\,\frac{\mu(y)}{\pi(y)}\right),\quad \{x,y\}\in E,
        \end{align*}
    with associated Laplacian matrix $L_{w^{(\mu)}}$. Let $h\in\mathbb{R}^{V}_0$ be fixed with $\|h\|_2\le 1$. Then there exist constants $C,C' > 0$, depending only on $G$ and $\mu$, such that for $\varepsilon\in\mathbb{R}$ with $|\varepsilon|$ sufficiently small (depending only on $G$ and $\mu$), it holds
        \begin{align*}
            \frac{\varepsilon^2}{1 + C'|\varepsilon|}\, h^\top L_{w^{(\mu)}}^{\dagger} h
            \le \Maas(\mu,\mu+\varepsilon h)^2
            \le \frac{\varepsilon^2}{1 - C|\varepsilon|}\, h^\top L_{w^{(\mu)}}^{\dagger} h.
        \end{align*}
    In particular,
        \begin{align*}
            \Maas(\mu,\mu+\varepsilon h)^2 &= |\varepsilon|^2 \, h^\top L_{w^{(\mu)}}^{\dagger} h + o(|\varepsilon|^2).
        \end{align*}
\end{theorem}

We prove~\Cref{thm:linearization} as~\Cref{thm:linearization-mu} in~\Cref{sec:linearization}. It is not difficult to obtain an asymptotic version of~\Cref{thm:linearization} using known facts about $\Maas(\cdot,\cdot)$ and its metric tensor, indeed, results reminiscent of this are already known in the literature (see, e.g.,~\cite{li2022transport}; we also discuss this in greater detail at the beginning of~\Cref{sec:linearization}). Thus our contribution is primarily the development of a \emph{nonasymptotic} version with explicit constants, which is arguably more suitable for applications. We state below a straightforward corollary in the simplified case of a stationary reference point.

\begin{corollary}[Stationary reference measure]
    Consider the setting described in~\Cref{thm:linearization}. If $\mu=\pi$, then $w^{(\pi)}=w$ and hence $L_{w^{(\pi)}}=L_w$. In this case,~\Cref{thm:linearization} specializes to the stationary linearization
        \begin{align}\label{eq:linearization-pi}
            \Maas(\pi,\pi+\varepsilon h)^2 = |\varepsilon|^2 \, h^\top L_{w}^{\dagger} h + o(|\varepsilon|^2),
        \end{align}
    with constants depending only on $G$.
\end{corollary}

We remark that the leading-order term in~\cref{eq:linearization-pi} can be considered a graph analogue of the dual homogeneous Sobolev norm, in which case~\cref{eq:linearization-pi} parallels the continuum $\dot H^{-1}$ linearization in~\cref{eq:sobolev-linearization-asymptotic} (see~\Cref{thm:sobolev-characterization} and the surrounding section for a more detailed discussion of this viewpoint). 

At a high level,~\Cref{thm:linearization} asserts that, in a neighborhood of a reference point $\mu\in\Prob_\star(V)$, the tangent space to $\Prob_\star(V)$ can be identified with $\mathbb{R}^{V}_0$ and that locally the metric tensor is given by the Hilbertian inner product
    \begin{align}\label{eq:inner-product-ldagger}
        \langle f, g\rangle_{L_{w^{(\mu)}}^{\dagger}} &:= f^\top L_{w^{(\mu)}}^{\dagger} g,\quad f, g\in\mathbb{R}^{V}_0.
    \end{align}
In particular, in the stationary case $\mu=\pi$ one has $w^{(\pi)}=w$ and hence $\langle\cdot,\cdot\rangle_{L_{w^{(\pi)}}^{\dagger}}=\langle\cdot,\cdot\rangle_{L_w^{\dagger}}$. The induced norm distance $\|f-g\|_{L_{w}^{\dagger}}^2 = (f-g)^\top L_{w}^{\dagger} (f-g)$ is actually well known when evaluated on Dirac masses: for $x,y\in V$, their \emph{effective resistance}, or simply resistance distance, is given by the formula
    \begin{align}\label{eq:resistance-distance-vertices}
        \Resis (x, y) &= (\mathbf{1}_x-\mathbf{1}_y)^\top L_w^\dagger (\mathbf{1}_x-\mathbf{1}_y).
    \end{align}
Here and throughout the paper, $\mathbf{1}_x\in\mathbb{R}^{V}$ denotes the indicator vector of the vertex $x\in V$, while $\mathbf{1}$ denotes the all ones vector. Thus~\cref{eq:inner-product-ldagger} can be seen as a straightforward extension of resistance distance from \emph{vertices} to general probability densities. An extension of this sort appears to have been first mentioned in passing in~\cite{jorgensen2010hilbert} but to the best of these authors' knowledge has not attracted deeper scrutiny in the literature.

\begin{definition}[Resistance distance between measures]\label{defn:resistance-distance-measures}
    Let $G=(V,E,w)$ be a weighted graph. For $\mu,\nu\in\Prob(V)$, we define
        \begin{align*}
            \Resis (\mu,\nu) := (\mu-\nu)^\top L_{w}^\dagger (\mu-\nu).
        \end{align*}
\end{definition}

We remark that in~\Cref{defn:resistance-distance-measures}, the term ``distance'' is used with slight abuse. It is well known that both $(x,y)\mapsto \Resis(x,y)$ and $(x,y)\mapsto \Resis(x,y)^{1/2}$ define metrics on $V$. However, $(\mu,\nu)\mapsto \Resis (\mu,\nu)$ fails to satisfy the triangle inequality in general, and thus $\Resis (\mu,\nu)$ defines only a semimetric on $\Prob(V)$ (it is nonnegative, symmetric, and vanishes if and only if $\mu=\nu$). We remark that the weaker bound $\Resis (\mu,\nu)\le 2\Resis (\mu,\tau) + 2\Resis (\tau,\nu)$ does hold for $\mu,\nu,\tau\in\Prob(V)$ and that, as in the vertex case, $(\mu,\nu)\mapsto \Resis(\mu,\nu)^{1/2}$ defines a true metric on $\Prob(V)$. 

We thusly investigate what we term the \emph{resistance manifold}, which is the metric space $(\Prob_\star(V),d_{\Resis})$, where $d_\Resis(\mu,\nu) := \Resis(\mu,\nu)^{1/2}$ for each $\mu,\nu\in\Prob_\star(V)$. 

\paragraph{Our second contribution.}
Upon viewing $(\Prob_\star(V),d_{\Resis})$ as a flat Riemannian manifold, our second main contribution is to characterize the resulting geodesic distance in a number of ways.  Notably, effective resistance is related to a number of \emph{combinatorial} properties of the underlying graph, including the spanning trees and $2$-forests of $G$~\cite{chaiken1982combinatorial,barrett2020spanning}, the eigenvalues of $L_{w}$ and its normalized counterpart~\cite{chen2007resistance}, and hitting times associated with the simple random walk on $V$~\cite{tetali1991random,lovasz1995mixing}. Thus the geodesic distance $d_{\Resis}$ can be viewed through this lens as well.

We define the \emph{volume} of $G$ as the sum of its vertex degrees: $\Vol_{w}(G) = \sum_{x\in V} \Deg_{x}$. If $F\subseteq E$, we define the weight of $F$, denoted $w_F$, by 
    \begin{align*}
        w_F &:= \prod_{e\in F} w_e.
    \end{align*}

\begin{theorem}[Characterizations of resistance distance between measures]\label{thm:characterizations}
    Let $G=(V,E, w)$ be a connected graph and let $\mu,\nu\in\Prob(V)$ be fixed. Then we have the following characterizations of $\Resis (\mu,\nu)$:
        \begin{enumerate}
            \item \emph{(Beckmann formula, see~\Cref{thm:beckmann})} Let $E'$ denote an orientation of $E$, and let $B\in\mathbb{R}^{V\times E'}$ denote the corresponding oriented incidence matrix. Then
                \begin{align*}
                    \Resis (\mu,\nu) = \inf\left\{\sum_{e=(x,y)\in E'} \frac{|J(e)|^2}{w_e} \;:\; J\in\mathbb{R}^{E'},\, BJ = \mu-\nu\right\}.
                \end{align*}

            \item \emph{(Benamou--Brenier formula, see~\Cref{thm:benamou-brenier-graphs})} Let $E'$ denote an orientation of $E$. Then
                \begin{align*}
                    \Resis (\mu,\nu) &= \inf_{\substack{\rho_t\in \mathcal{C}^1([0,1]; \mathbb{R}^{V})\\ J_t\in \mathcal{B}([0,1]; \mathbb{R}^{E'})}} \left\{\int_0^1 \sum_{e\in E'} \frac{|J_t(e)|^2}{w_e} \, dt \,:\, \frac{d}{dt} \rho_t + B J_t = 0,\quad \rho_0 = \mu,\quad \rho_1 = \nu\right\}.
                \end{align*}

            \item \emph{(Negative homogeneous Sobolev distance, see~\Cref{thm:sobolev-characterization})}  $\Resis (\mu,\nu) = \|\mu-\nu\|_{\dot{H}^{-1}(V,w)}^2$, where $\|\cdot\|_{\dot{H}^{-1}(V,w)}$ is the norm dual to
                \begin{align*}
                    \|f\|_{\dot{H}^{1}(V,w)}^{2} &:= \sum_{\{x,y\}\in E} w_{xy} |f(x)-f(y)|^2 ,\quad f\in \mathbb{R}^{V}_0.
                \end{align*}
    
            \item \emph{(Spanning $2$-forest formula, see~\Cref{thm:measure-spanning-two-forest})} 
                \begin{align}\label{eq:two-forest-characterization}
                    \Resis (\mu,\nu) &= \frac{1}{\Tree (G)}\sum_{\substack{F\in\Forest (G) \\ F=T_1\cup T_2}} w_F \, (\mu(T_1)-\nu(T_1))^2,
                \end{align}
            where $\Tree(G)$ denotes the sum of all weights $w_T$ of spanning trees $T$ of $G$, and the sum on the right-hand side of~\cref{eq:two-forest-characterization} runs over the set of spanning $2$-forests $F$ of $G$, denoted $\Forest (G)$, with components $T_1$ and $T_2$.

            \item \emph{(Hitting time formula, see~\Cref{thm:resistance-access-time})} For $x, y\in V$, let $\Hitting(x, y)$ denote the expected number of steps for the simple random walk on $V$ to reach $y$ having started at $x$. Then
                \begin{align*}
                    \Resis (\mu,\nu) &= -\frac{1}{\Vol_w(G)} \sum_{x, y\in V}(\mu(x) - \nu(x)) (\mu(y) - \nu(y)) \Hitting(x, y).
                \end{align*}
        \end{enumerate}
\end{theorem}

\Cref{thm:characterizations} is proved in~\Cref{sec:characterizations}. 

We explore more connections between the combinatorial structure of the underlying graph $G$, the discrete transportation metric $\Maas(\cdot,\cdot)$, and its linearization by considering the concrete setting of bipartite graphs in~\Cref{subsec:bipartite}. In this environment we can provide quantitative descriptions of both $\Maas(\cdot,\cdot)$ and $\Resis(\cdot,\cdot)$ for a one-parameter family of perturbations of the stationary distribution of the simple random walk.

\paragraph{Our third contribution.} We further investigate the geometry of the resistance manifold through the lens of gradient flows. Inspired by Maas' result that entropy gradient flow on the probability simplex recovers the diffusion generated by the underlying random walk, we show that the $\chi^{2}$ functional on the flat manifold $(\Prob_\star(V),d_\Resis)$ has gradient flow given by the continuous time random walk. Relatedly, we are able to control the relaxation rate of this gradient flow by the eigenvalues of the normalized Laplacian matrix of $G$, thereby recovering a classical result from spectral graph theory. Taken together, these observations constitute our third contribution, stated in full below.

If $G$ is connected and $\pi\in\Prob(V)$ denotes the stationary distribution of the simple random walk on $V$, then we define the (normalized) \emph{$\chi^2$ functional} on $\Prob(V)$ according to the formula
    \begin{align}\label{eq:chi2-functional-sec4}
        \chi^2(\mu) &= \frac{1}{2\,\Vol_w(G)}\sum_{x\in V} \frac{(\mu(x) - \pi(x))^2}{\pi(x)},\quad \mu\in\Prob(V).
    \end{align} 

\begin{theorem}[Gradient flows of the $\chi^{2}$ functional]\label{thm:gradient-flow}
    Let $G=(V,E,w)$ be a connected graph and let $\pi\in\Prob(V)$ denote the stationary distribution of the simple random walk on $V$. Then, on $\Prob_\star(V)$, the gradient flow of $\chi^2$ with respect to the Riemannian metric tensor induced by the inner product $\langle\cdot,\cdot\rangle_{L_{w}^\dagger}$ coincides with the continuous time evolution of the random walk.
\end{theorem}

\Cref{thm:gradient-flow} is proved in~\Cref{sec:gradient-flows} as~\Cref{cor:chi2-fdiv}.

To state our final result, we recall the notion of geodesic strong convexity. If $(\mathcal{X}, d)$ is a Riemannian manifold, a $\mathcal{C}^2$ functional $F:\mathcal{X}\to\mathbb{R}$ is \emph{geodesic $\lambda$-strongly convex} if, for every pair $\rho_0,\rho_1\in\mathcal{X}$ and every constant-speed geodesic $(\rho_t)_{t\in[0,1]}$ joining them, one has:
    \begin{align}\label{eq:geodesic-convexity}
        F(\rho_t) &\leq (1-t)F(\rho_0) + t F(\rho_1) - \frac{\lambda}{2} t(1-t) d(\rho_0,\rho_1)^2
    \end{align}
for each $0\le t\le 1$. The greatest such $\lambda\in\mathbb{R}$ for which~\cref{eq:geodesic-convexity} holds is known as the \emph{geodesic strong convexity modulus} of the functional $F$. 

We denote by $\mathcal{L}_w = D_w^{-1/2} L_w D_w^{-1/2}$ the \emph{normalized Laplacian matrix} of $G$, and recall that $\mathcal{L}_w$ is symmetric positive semidefinite with nullspace $\mathbb{R} D_w^{1/2}\mathbf{1}$. The least strictly positive eigenvalue of $\mathcal{L}_w$ is denoted $\lambda_2(\mathcal{L}_w)$ and is called the \emph{spectral gap} of $\mathcal{L}_w$.

\begin{theorem}[Geodesic strong convexity of the $\chi^{2}$ functional]\label{thm:modulus}
    Let $G=(V,E,w)$ be a connected graph and let $\pi$ denote the stationary distribution of the associated simple random walk on $V$. On the Riemannian manifold $(\Prob_\star(V),d_{\Resis})$, the functional $\chi^2$ has geodesic strong convexity modulus equal to the spectral gap $\lambda_2(\mathcal{L}_w)$ of $\mathcal{L}_{w}$.
\end{theorem}

\Cref{thm:modulus} is proved in~\Cref{sec:gradient-flows} as~\Cref{thm:modulus-sec4}.

We remark that geodesic $\lambda$-strong convexity yields exponential relaxation of the associated gradient flow, with $d(\rho_t,\rho_\star)\lesssim e^{-\lambda t}$ and $|F(\rho_t)-F(\rho_\star)|\lesssim e^{-2\lambda t}$ for a minimizer $\rho_\star$ (see, e.g.,~\cite[Sec. 11.2]{ambrosio2005gradient}). In our setting,~\Cref{thm:gradient-flow} and~\Cref{thm:modulus} therefore recover the classical result concerning $\chi^2$-relaxation of the random walk diffusion process at a rate governed by the spectral gap $\lambda_2(\mathcal{L}_w)$ (see, e.g.,~\cite{chung1997spectral}).

From a computational perspective, our results suggest that the discrete transportation metric may be well approximated in certain regimes by the surrogate $\Resis(\cdot,\cdot)$ whose computation amounts to basic linear algebra and which is substantially more tractable than solving the full dynamic transport problem. While we do not develop algorithms directly here, our results develop a geometric framework that naturally supports scalable approximations, local embeddings, and various numerical schemes for graph transport. We feel these topics comprise promising directions for future investigation.

\subsection{Related work}\label{subsec:related-work}

Research concerning optimal transport on graphs can be modeled as having developed along two related but distinct directions. One such thread considers classical Wasserstein transportation problems using the ambient graph metric, and includes topics such as explicit formulas on trees, algorithmic approximations and regularization, and scalable computation on mesh discretizations \cite{evans2012phylogenetic,mcgregor2013sketching,solomon2014earth,essid2018quadratically}. Other directions include Sobolev transport on graphs~\cite{le2022sobolev,le2023scalable,le2024generalized} as well as transport for vector-valued densities~\cite{robertson2025generalization,chen2018vector}. We also remark that $\Wass_{1}$ on graphs is related to (uncapacitated) minimum cost flow programming (see, e.g.,~\cite[Ch. 6]{peyre2019computational}), for which there is an extensive literature that we will not summarize here.

Another direction (see~\cite{maas2011gradient,mielke2011gradient,li2022transport}), formulates transportation using action functionals of the form~\cref{eq:action-functional}. Subsequent work has included the analysis of gradient flows for Markov chains~\cite{chow2012fokker}, lazy random walks and Schr\"odinger-type dynamics~\cite{leonard2015lazy,chow2019discrete}, and discrete Ricci curvature and entropy convexity theory~\cite{erbar2012ricci,mielke2013geodesic,erbar2014gradient}. It has also motivated numerical and asymptotic analyses, including computational results on finite graphs, convergence to continuum Wasserstein distances, and extensions to metric graphs~\cite{erbar2020computation,gigli2013gromov,erbar2022gradient}.

Meanwhile, effective resistance dates back roughly to the 1980s and 1990s when it was developed as a tool for analyzing random walks and electrical networks on graphs; see, e.g.,~\cite{doylesnell1984random,chandra1996electrical,tetali1991random}. It has since emerged as a fundamental object in modern computational graph theory with an extensive literature touching on a number of areas; while a full survey is outside the scope of the present paper, we highlight well-known applications in clustering and recommendation~\cite{fouss2007random}, network statistics~\cite{newman2005measure,brandes2005centrality}, and network design and spectral sparsification~\cite{ghosh2008minimizing,spielman2011graph}. We also remark that effective resistance has been investigated in settings beyond weighted graphs, including directed graphs~\cite{young2016newi,young2016newii,sugiyama2023kron} and graphs with matrix-valued weights~\cite{atik2019resistance,cloninger2024random}.

A substantial part of the linearized optimal transport literature concerns embeddings built from Monge maps relative to a fixed reference measure. In this line, Wang et al.\ and Kolouri et al.\ introduced linear transport embeddings for image analysis and pattern recognition by representing each measure through its transport to a reference image \cite{wang2013linear,kolouri2016continuous}. Subsequent work kept this reference map viewpoint and elucidated the circumstances under which it behaves as a controlled surrogate for $\Wass_2$: M\'erigot, Delalande, and Chazal proved quantitative stability of the embedding, while Delalande and M\'erigot sharpened this analysis under perturbations of the target \cite{merigot2020quantitative,delalanda2023quantitative}. Additional work has developed the same general paradigm for classification guarantees, statistical rates, singular reference measures, manifolds, and other transport geometries \cite{park2018cumulative, moosmuller2023linear,aldroubi2021partitioning, cloninger2025linearized, nenna2023transport,sarrazin2024linearized,cai2022linearized,bai2023linear}. 

A distinct line of work, closer in spirit to the present paper, does not encode measures by Monge maps to a fixed base. Instead it works directly with the quadratic form on the tangent space, typically expressed through elliptic operators and negative homogeneous Sobolev norms. In the continuum, the comparison results of Peyr\'e and the weighted linearization investigated in a preprint of Greengard et al.\ make this relation explicit \cite{peyre2018comparison,greengard2022linearization}. Related Sobolev models have also been used in practice in inverse data matching and seismic inversion \cite{engquist2020quadratic,engquist2022optimal}, image denoising~\cite{zhu2024implicit}, and learning submanifolds of the Wasserstein manifold \cite{hamm2025manifold}. Our construction can be viewed as belonging to this line of work; to our knowledge, ours is the first work to provide quantitative guarantees for linearized transport on graphs. 

\subsection{Organization}

This paper is organized as follows. In~\Cref{sec:linearization} we develop the linearization of the discrete transportation metric, from the asymptotic viewpoint (\Cref{subsec:linearization-asymptotic}) to nonasymptotic bounds with explicit constants (\Cref{subsec:linearization-nonasymptotic}), and finally in the example setting of bipartite graphs (\Cref{subsec:bipartite}). In~\Cref{sec:characterizations} we prove the variational and combinatorial characterizations of $\Resis(\cdot,\cdot)$ stated in \Cref{thm:characterizations} (\Cref{subsec:variational-characterizations,subsec:combinatorial-characterizations}). Finally, in~\Cref{sec:gradient-flows} we study gradient flows on the resistance manifold, focusing specifically on $f$-divergences, and prove~\Cref{thm:gradient-flow} and~\Cref{thm:modulus}.

\section{Linearization of the discrete transportation distance}\label{sec:linearization}

This section focuses on our main result~\Cref{thm:linearization}. In~\Cref{subsec:linearization-asymptotic}, we detail the geometric setting of the result and discuss it in asymptotic terms. Then in~\Cref{subsec:linearization-nonasymptotic} we state and prove the nonasymptotic version. Finally, in~\Cref{subsec:bipartite}, as an example we explore the setting of bipartite graphs where, for a specific family of measures, the resistance distance and discrete transportation distance can be computed exactly and compared quantitatively.

\subsection{Linearization of the transportation metric in asymptotic terms}\label{subsec:linearization-asymptotic}

Letting $G=(V,E,w)$ be a fixed connected graph, we recall the definition of the probability simplex on $V$:
    \begin{align*}
        \Prob(V) &:= \left\{ \mu \in\mathbb{R}^{V} \,:\, \mu(x)\geq 0\text{ for each }x\in V,\, \sum_{x\in V} \mu(x) = 1\right\}.
    \end{align*}
Once again we write $\Prob_\star(V)$ to denote the interior of $\Prob(V)$. Letting $\mu\in\Prob_\star(V)$ be a fixed reference measure, it is straightforward to see that the embedded tangent space to $\Prob(V)$ at $\mu$ is given by a copy of the hyperplane $\mathbb{R}^{V}_0$. On the other hand, Maas~\cite{maas2011gradient} shows that if $\mu\in\Prob_\star(V)$, then by defining the edge weights
    \begin{align*}
        w^{(\mu)}_{xy} := w_{xy}\,\theta\left(\frac{\mu(x)}{\pi(x)},\,\frac{\mu(y)}{\pi(y)}\right),\quad \{x,y\}\in E,
    \end{align*}
and the associated Laplacian matrix $L_{w^{(\mu)}}$, then the intrinsic tangent space to $(\Prob_\star(V),\Maas)$ at $\mu$ can be identified with the set of discrete gradients
    \begin{align*}
        \{\nabla \psi \in\mathbb{R}^{V\times V} \,:\, \psi\in\mathbb{R}^{V},\, \nabla \psi(x,y) = \psi(y) - \psi(x)\}.
    \end{align*}
It is straightforward to verify that the map $\nabla \psi \mapsto L_{w^{(\mu)}} \psi$ defines a linear isomorphism from the intrinsic tangent space of $(\Prob_\star(V),\Maas)$ at $\mu$ to the embedded tangent space $\mathbb{R}^{V}_0$. Moreover, since the Riemannian metric on $(\Prob_\star(V),\Maas)$ at $\mu$ can be expressed in terms of $L_{w^{(\mu)}}$ via the formula
    \begin{align*}
        \langle \nabla \psi, \nabla \phi\rangle_\mu &= \psi^\top L_{w^{(\mu)}} \phi,
    \end{align*}
it follows that the linearization of the discrete transportation distance at $\mu$ is given by the inner product $\langle\cdot,\cdot\rangle_{L_{w^{(\mu)}}^\dagger}$ on $\mathbb{R}^{V}_0$ as defined in~\cref{eq:inner-product-ldagger}. Therefore if $\rho \in \mathcal{C}^2((-\varepsilon,\varepsilon);\Prob_\star(V))$ is any curve satisfying $\rho_0 = \mu$ and $\frac{d}{dt}\rho|_{t=0}=:h$, then since $t\mapsto \Maas(\mu,\rho_t)^2$ is $\mathcal{C}^{2}$ in a neighborhood of $t=0$ (see, e.g.,~\cite[Lemma 6.8]{lee2018introduction}), it follows from Taylor's theorem that 
    \begin{align}\label{eq:asymptotic-linearization}
        \Maas(\mu,\rho_t)^2 &= t^2 \, h^\top L_{w^{(\mu)}}^\dagger h + o(t^2),\quad t\to 0.
    \end{align}

\subsection{Linearization of the transportation metric in nonasymptotic terms}\label{subsec:linearization-nonasymptotic}

\begin{theorem}\label{thm:linearization-mu}
    Let $G=(V,E, w)$ be a connected graph and let $\pi\in\Prob_\star (V)$ denote the stationary distribution of the simple random walk on $V$, and let $\mu\in\Prob_\star(V)$ be fixed. Define edge weights $w^{(\mu)}=(w^{(\mu)}_{xy})_{\{x,y\}\in E}$ by
        \begin{align*}
            w^{(\mu)}_{xy} := w_{xy}\,\theta\left(\frac{\mu(x)}{\pi(x)},\,\frac{\mu(y)}{\pi(y)}\right),\quad \{x,y\}\in E,
        \end{align*}
    and let $L_{w^{(\mu)}}$ be the corresponding Laplacian. Let $h\in\mathbb{R}^{V}_0$ be fixed with $\|h\|_2\le 1$. Let $\lambda_2(L_{w^{(\mu)}})$ denote the smallest nonzero eigenvalue of $L_{w^{(\mu)}}$, and define constants $0 < C \le C'$ via
        \begin{align}\label{eq:constants-linearization}
            C &:= \left(\min_{x\in V}\mu(x)\right)^{-1},\quad C' := C\,\sqrt{\frac{2\,\Vol_w(G)}{\lambda_2(L_{w^{(\mu)}})}}.
        \end{align}
    Then for any $\varepsilon\in\mathbb{R}$ with $|\varepsilon|$ sufficiently small (depending only on $G$ and $\mu$), it holds
        \begin{align}\label{eq:main-linearization}
            \frac{\varepsilon^2}{1 + C'|\varepsilon|}\, h^\top L_{w^{(\mu)}}^\dagger h
            \le \Maas(\mu,\mu+\varepsilon h)^2
            \le \frac{\varepsilon^2}{1 - C|\varepsilon|}\, h^\top L_{w^{(\mu)}}^\dagger h.
        \end{align}
    In particular, $\Maas(\mu,\mu+\varepsilon h)^2 = |\varepsilon|^2\, h^\top L_{w^{(\mu)}}^\dagger h  + o(|\varepsilon|^2)$.
\end{theorem}

\begin{proof}
    First we will exhibit an admissible pair $(\rho_t), (\psi_t)$ as in~\cref{eq:defn-maas-distance} to obtain the upper bound in~\cref{eq:main-linearization}. To this end let
        \begin{align*}
            \rho_t &= t\, (\mu + \varepsilon h) + (1-t) \, \mu,\quad 0\le t\le 1.
        \end{align*}
    Then if $|\varepsilon|$ is chosen to be small enough, $(\rho_t)\in \mathcal{C}^{1}([0, 1];\,\Prob_\star(V))$. Clearly $\rho_0 = \mu,\, \rho_1=\mu + \varepsilon h$. The flux constraint in~\cref{eq:defn-maas-distance} requires of $\psi_t$ that
        \begin{align}\label{eq:flux-constraint-linear}
            \frac{d}{dt} \rho_t(x) &= \varepsilon h(x) = -\sum_{y\sim x} w_{xy}\, (\psi_t(y) - \psi_t(x))\, \theta\left(\frac{\rho_t(x)}{\pi(x)}, \frac{\rho_t(y)}{\pi(y)}\right),\quad x\in V,\, 0\le t\le 1.
        \end{align} 
    Write
        \begin{align*}
            \widetilde{w}^{(t)}_{xy} &= w_{xy}\,\theta\left(\frac{\rho_t(x)}{\pi(x)}, \frac{\rho_t(y)}{\pi(y)}\right),\quad \{x, y\}\in E,
        \end{align*}
    so that~\cref{eq:flux-constraint-linear} becomes
        \begin{align}\label{eq:flux-constraint-linear-2}
            \varepsilon \, h &= L_{\widetilde{w}^{(t)}} \psi_t,\quad 0\le t\le 1.
        \end{align}
    Note that since $\rho_t\in\Prob_\star (V)$ holds for all $t$ and since $\theta(a, b) > 0$ whenever $a, b > 0$, the weights $\widetilde{w}^{(t)}_{xy}$ remain strictly positive for all $0\le t\le 1$. In particular, the operator $L_{\widetilde{w}^{(t)}}$ remains the Laplacian of a connected weighted graph and is therefore positive definite on the subspace $\mathbb{R}^{V}_0$ which includes $h$. We may therefore pick, for each $0\le t\le 1$, a unique representative $\psi_{t}^\ast\in\mathbb{R}^{V}_0$ satisfying~\cref{eq:flux-constraint-linear-2} and given by
        \begin{align*}
            \psi_t^\ast &= \varepsilon \, L_{\widetilde{w}^{(t)}}^{\dagger} h.
        \end{align*}
    Then, using $L_{\widetilde{w}^{(t)}}^{\dagger}L_{\widetilde{w}^{(t)}} L_{\widetilde{w}^{(t)}}^{\dagger} = L_{\widetilde{w}^{(t)}}^{\dagger}$, the action functional evaluates to
        \begin{align}\label{eq:action-upper-bound-1}
            \mathcal{A}((\psi_t^\ast), (\rho_t)) &= \int_0^{1} (\psi_t^\ast)^\top L_{\widetilde{w}^{(t)}} \psi_t^\ast \, dt = \varepsilon^2 \int_0^{1} h^\top L_{\widetilde{w}^{(t)}}^{\dagger} h \, dt\notag\\
            &= \varepsilon^2 \, h^\top \left( \int_0^{1} L_{\widetilde{w}^{(t)}}^{\dagger} \, dt \right) h.
        \end{align}
    To bound~\cref{eq:action-upper-bound-1} from above we first observe that, for each $x\in V$ and $0\le t\le 1$, one has
        \begin{align*}
            \frac{\rho_t(x)}{\pi(x)} &= \frac{\mu(x) + t \varepsilon h(x)}{\pi(x)} = \frac{\mu(x)}{\pi(x)} + t \varepsilon \frac{h(x)}{\pi(x)}.
        \end{align*}
    Using $|h(x)|\le \|h\|_2\le 1$, we have
        \begin{align*}
            \left|t \varepsilon \frac{h(x)}{\pi(x)}\,\frac{\pi(x)}{\mu(x)}\right| &\leq |\varepsilon| C,
        \end{align*}
    where $C$ is taken as in~\cref{eq:constants-linearization}. Thus it follows that
        \begin{align}\label{eq:pi-mu-bound}
            \left(1-|\varepsilon|C\right)\frac{\mu(x)}{\pi(x)} \le \frac{\rho_t(x)}{\pi(x)} \le \left(1+|\varepsilon|C\right)\frac{\mu(x)}{\pi(x)},\quad x\in V,\,0\le t\le 1.
        \end{align}
    Since $\theta$ is increasing in each argument and homogeneous of degree one,~\cref{eq:pi-mu-bound} yields
        \begin{align*}
            \left(1-|\varepsilon|C\right)\theta\left(\frac{\mu(x)}{\pi(x)},\frac{\mu(y)}{\pi(y)}\right)
            \le \theta\left(\frac{\rho_t(x)}{\pi(x)},\frac{\rho_t(y)}{\pi(y)}\right)
            \le \left(1+|\varepsilon|C\right)\theta\left(\frac{\mu(x)}{\pi(x)},\frac{\mu(y)}{\pi(y)}\right),
        \end{align*}
    and therefore, for each $\{x, y\}\in E$ and $0\le t\le 1$, it holds
        \begin{align*}
            (1 - |\varepsilon| C) \, w^{(\mu)}_{xy} \leq \widetilde{w}^{(t)}_{xy} \leq (1 + |\varepsilon| C) \, w^{(\mu)}_{xy}.
        \end{align*}
    Now let $f\in\mathbb{R}^{V}$ be arbitrary. It follows that
        \begin{align*}
            (1 - |\varepsilon| C) \, f^\top L_{w^{(\mu)}} f &\leq f^\top L_{\widetilde{w}^{(t)}} f \leq (1 + |\varepsilon| C) \, f^\top L_{w^{(\mu)}} f,
        \end{align*}
    and in particular this bound holds for $f\in\mathbb{R}^{V}_0$. Since $L_{w^{(\mu)}}$ and $L_{\widetilde{w}^{(t)}}$ are positive definite on $\mathbb{R}^{V}_0$, the inverse inequality holds for their Moore--Penrose inverses:
        \begin{align}\label{eq:pseudoinverse-bound}
            \frac{1}{1 + |\varepsilon| C} \, f^\top L_{w^{(\mu)}}^{\dagger} f &\leq f^\top L_{\widetilde{w}^{(t)}}^{\dagger} f \leq \frac{1}{1 - |\varepsilon| C} \, f^\top L_{w^{(\mu)}}^{\dagger} f,\quad f\in\mathbb{R}^{V}_0.
        \end{align}
    Since~\cref{eq:pseudoinverse-bound} holds for all $0\le t\le 1$, we may integrate all sides to obtain
        \begin{align*}
            \frac{1}{1 + |\varepsilon| C} \, f^\top L_{w^{(\mu)}}^{\dagger} f &\leq f^\top \left( \int_0^{1} L_{\widetilde{w}^{(t)}}^{\dagger} \, dt \right) f \leq \frac{1}{1 - |\varepsilon| C} \, f^\top L_{w^{(\mu)}}^{\dagger} f,
        \end{align*}
    thus, with $f=h$ and via~\cref{eq:action-upper-bound-1}, one has
        \begin{align}\label{eq:pseudoinverse-bound-2}
            \mathcal{A}((\psi_t^\ast), (\rho_t)) &= \varepsilon^2 \, h^\top \left( \int_0^{1} L_{\widetilde{w}^{(t)}}^{\dagger} \, dt \right) h \leq \frac{\varepsilon^2}{1-|\varepsilon| C} h^\top L_{w^{(\mu)}}^{\dagger} h.
        \end{align}
    
    Next we prove the lower bound. By~\cref{eq:pseudoinverse-bound-2} and taking $|\varepsilon|C\le \frac12$, we have
        \begin{align*}
            \Maas(\mu, \mu+\varepsilon \, h)^2 &\leq 2\varepsilon^2 h^\top L_{w^{(\mu)}}^{\dagger} h 
        \end{align*}
    Since $\|h\|_2\le 1$ and $h\in\mathbb{R}^{V}_0$, we have $h^\top L_{w^{(\mu)}}^{\dagger} h \le \lambda_2(L_{w^{(\mu)}})^{-1}\|h\|_2^2 \le \lambda_2(L_{w^{(\mu)}})^{-1}$, so
        \begin{align*}
            \Maas(\mu, \mu+\varepsilon \, h)^2 &\leq \frac{2}{\lambda_2(L_{w^{(\mu)}})}\varepsilon^2.
        \end{align*}
    Now let $\eta > 0$ be fixed, and let $(\rho_t), (\psi_t)$ be any admissible pair of curves satisfying
        \begin{align*}
            \mathcal{A}((\psi_t), (\rho_t)) &\leq \Maas(\mu, \mu+\varepsilon h)^2 + \eta \leq \frac{2}{\lambda_2(L_{w^{(\mu)}})}\varepsilon^2 + \eta.
        \end{align*}
    We emphasize that $(\rho_t), (\psi_t)$ depend on the choice of $\eta$, but we suppress this in the notation. Then for any $x\in V$ and $0\le t\le 1$, since $(\rho_t)$ is piecewise $\mathcal{C}^1$, one has that
        \begin{align*}
            \rho_t(x) - \rho_0(x) &= \int_0^t \frac{d}{ds} \rho_s(x) \, ds = \int_0^t L_{\widehat{w}^{(s)}} \psi_s(x) \, ds,
        \end{align*}
    where $L_{\widehat{w}^{(s)}}$ is the Laplacian with edge weights
        \begin{align*}
            \widehat{w}^{(s)}_{xy} = w_{xy}\, \theta\left(\frac{\rho_s(x)}{\pi(x)}, \,\frac{\rho_s(y)}{\pi(y)}\right).
        \end{align*}
    By the Cauchy--Schwarz inequality, with $\rho_0 = \mu$,
        \begin{align}\label{eq:cauchy-schwarz-1}
            |\rho_t(x) - \mu(x)|^2 &\leq \left|\int_0^t L_{\widehat{w}^{(s)}} \psi_s(x) \, ds \right|^2 \leq t \int_0^t |L_{\widehat{w}^{(s)}} \psi_s(x)|^2 \, ds \leq \int_0^1 |L_{\widehat{w}^{(s)}} \psi_s(x)|^2 \, ds.
        \end{align}
    Let $x\in V$ be fixed. Then, once again by Cauchy--Schwarz and the definition of $\widehat{w}^{(s)}_{xy}$, we have
        \begin{align*}
            |L_{\widehat{w}^{(s)}} \psi_s(x)|^2 &= \left|\sum_{y\sim x} \widehat{w}^{(s)}_{xy} (\psi_s(y) - \psi_s(x))\right|^2 \\
            &\leq  \left(\sum_{y\sim x} \widehat{w}^{(s)}_{xy} \right)\left(\sum_{y\sim x} \widehat{w}^{(s)}_{xy} (\psi_s(y) - \psi_s(x))^2 \right)\\
            &= \left(\sum_{y\sim x} w_{xy}\, \theta\left(\frac{\rho_s(x)}{\pi(x)}, \frac{\rho_s(y)}{\pi(y)}\right) \right)\left(\sum_{y\sim x} \widehat{w}^{(s)}_{xy} (\psi_s(y) - \psi_s(x))^2 \right).
        \end{align*}
    Define $r_z := \rho_s(z)/\pi(z)$ for $z\in V$. Using the bound $\theta(a,b)\le (a+b)/2$ for $a,b\ge 0$, we obtain
        \begin{align*}
            \sum_{y\sim x}\widehat{w}^{(s)}_{xy}
            &= \sum_{y\sim x} w_{xy}\,\theta(r_x,r_y)
            \le \frac{1}{2}r_x\sum_{y\sim x}w_{xy} + \frac{1}{2}\sum_{y\sim x}w_{xy}r_y.
        \end{align*}
    Since $\pi(x)=\Deg_x/\Vol_w(G)$, the first term satisfies
        \begin{align*}
            \frac{1}{2}r_x\sum_{y\sim x}w_{xy} = \frac{1}{2}r_x\Deg_x=\frac{1}{2}\Vol_w(G)\rho_s(x) \le \frac{1}{2}\Vol_w(G).
        \end{align*}
    For the second term, using $w_{xy}\le \Deg_y$, we have
        \begin{align*}
            \sum_{y\sim x}w_{xy}r_y
            &= \Vol_w(G)\sum_{y\sim x}\frac{w_{xy}}{\Deg_y}\rho_s(y)
            \le \Vol_w(G)\sum_{y\in V}\rho_s(y)=\Vol_w(G).
        \end{align*}
    Therefore, for each $x\in V$ and $0\le s\le 1$,
        \begin{align*}
            \sum_{y\sim x}\widehat{w}^{(s)}_{xy} \le \Vol_w(G).
        \end{align*}
    Plugging this into the bound for $|L_{\widehat{w}^{(s)}} \psi_s(x)|^2$, we have
        \begin{align}\label{eq:estimate-for-intL}
            \int_0^1 |L_{\widehat{w}^{(s)}} \psi_s(x)|^2\,ds &\leq \Vol_w(G)\int_0^1 \sum_{y\sim x} \widehat{w}^{(s)}_{xy} (\psi_s(y) - \psi_s(x))^2\,ds\notag\\
            &\leq \Vol_w(G)\int_0^1 \sum_{\{u,v\}\in E} \widehat{w}^{(s)}_{uv} (\psi_s(u) - \psi_s(v))^2\,ds\notag\\
            &= \Vol_w(G)\mathcal{A}((\psi_t), (\rho_t)).
        \end{align}
    Recalling~\cref{eq:cauchy-schwarz-1} and using~\cref{eq:estimate-for-intL}, we obtain
        \begin{align*}
            |\rho_t(x) - \mu(x)|^2 &\leq \int_0^1 |L_{\widehat{w}^{(s)}} \psi_s(x)|^2 \, ds \leq \Vol_w(G) \mathcal{A}((\psi_t), (\rho_t)) \leq \Vol_w(G)\left(\frac{2}{\lambda_2(L_{w^{(\mu)}})}\varepsilon^2 + \eta\right),
        \end{align*}
    and in particular it follows that
        \begin{align}\label{eq:bound-rho-mu}
            \frac{\rho_t(x)}{\mu(x)} &\leq 1 + \frac{1}{\mu(x)}\sqrt{\Vol_w(G)\left(\frac{2}{\lambda_2(L_{w^{(\mu)}})}\varepsilon^2 + \eta\right)}\notag\\
            &\leq 1 + \left(\min_{u\in V}\mu(u)\right)^{-1}\sqrt{\Vol_w(G)}\left(\sqrt{\frac{2}{\lambda_2(L_{w^{(\mu)}})}}\,|\varepsilon| + \sqrt{\eta}\right).
        \end{align}
    Equivalently, defining
        \begin{align*}
            C' &:= \left(\min_{u\in V}\mu(u)\right)^{-1}\sqrt{\frac{2\,\Vol_w(G)}{\lambda_2(L_{w^{(\mu)}})}},\\
            C''_{\eta} &:= \left(\min_{u\in V}\mu(u)\right)^{-1}\sqrt{\Vol_w(G)}\sqrt{\eta},
        \end{align*}
    we may write~\cref{eq:bound-rho-mu} in the form
        \begin{align*}
            \frac{\rho_t(x)}{\mu(x)} &\leq 1 + C' |\varepsilon| + C''_{\eta}.
        \end{align*}
    Thus, for each $\{x,y\}\in E$, we have
        \begin{align*}
            \theta\left(\frac{\rho_t(x)}{\pi(x)}, \frac{\rho_t(y)}{\pi(y)}\right)
            &= \theta\left(\frac{\mu(x)}{\pi(x)}\frac{\rho_t(x)}{\mu(x)},\,\frac{\mu(y)}{\pi(y)}\frac{\rho_t(y)}{\mu(y)}\right)\\
            &\leq \theta\left(\left(1+C'|\varepsilon|+C''_{\eta}\right)\frac{\mu(x)}{\pi(x)},\,\left(1+C'|\varepsilon|+C''_{\eta}\right)\frac{\mu(y)}{\pi(y)}\right)\\
            &= \left(1+C'|\varepsilon|+C''_{\eta}\right)\theta\left(\frac{\mu(x)}{\pi(x)},\,\frac{\mu(y)}{\pi(y)}\right),
        \end{align*}
    where the last step uses the fact that $\theta$ is 1-homogeneous. We thus have the bound, for each $f\in\mathbb{R}^{V}$ and $0\le t\le 1$,
        \begin{align*}
            f^\top L_{\widehat{w}^{(t)}} f &\leq \left(1 + C'|\varepsilon| + C''_{\eta}\right) f^\top L_{w^{(\mu)}} f,
        \end{align*} 
    and inversely so for $f\in\mathbb{R}^{V}_0$:
        \begin{align}\label{eq:bound-Z}
            f^\top L_{\widehat{w}^{(t)}}^{\dagger} f &\geq \frac{1}{1 + C'|\varepsilon| + C''_{\eta}} f^\top L_{w^{(\mu)}}^{\dagger} f.
        \end{align}
    Let $0\le t\le 1$ be fixed. The flux constraint requires that
        \begin{align*}
            \frac{d}{dt}\rho_t &= L_{\widehat{w}^{(t)}} \psi_t,
        \end{align*}
    and we may freely assume that $\psi_t\in\mathbb{R}^{V}_0$ since adding a constant to $\psi_t$ does not change the action or violate the flux constraint. Thus, since $L_{\widehat{w}^{(t)}}$ is positive definite on $\mathbb{R}^{V}_0$, $\psi_t$ is uniquely determined via
        \begin{align}\label{eq:psi-t-uniquely-determined}
            \psi_t &= L_{\widehat{w}^{(t)}}^{\dagger} \rho_{t}',
        \end{align}
    where $\rho_{t}' = \frac{d}{dt} \rho_t$ for brevity. We then have that, since $ L_{\widehat{w}^{(t)}}^{\dagger}  L_{\widehat{w}^{(t)}}  L_{\widehat{w}^{(t)}}^{\dagger}  =  L_{\widehat{w}^{(t)}}^{\dagger} $ and since $\rho_{t}'\in\mathbb{R}^{V}_0$ for all $t$, it holds
        \begin{align*}
            \mathcal{A}((\psi_t), (\rho_t)) &= \int_0^1 \psi_t^\top L_{\widehat{w}^{(t)}} \psi_t \, dt = \int_0^1 \left(\rho_{t}'\right)^\top L_{\widehat{w}^{(t)}}^{\dagger} \left(\rho_{t}'\right) dt \geq \frac{1}{1 + C'|\varepsilon| + C''_{\eta}} \int_0^1 (\rho_t')^\top L_{w^{(\mu)}}^\dagger (\rho_t') dt,
        \end{align*}
    where the last step follows from~\cref{eq:bound-Z}. Separately, by the Cauchy--Schwarz inequality,
        \begin{align*}
            \left(\int_0^1 \rho_t' \,dt\right)^\top L_{w^{(\mu)}}^{\dagger} \left(\int_0^1 \rho_t' \,dt\right) &= \left\|\int_0^1 (L_{w^{(\mu)}}^{\dagger})^{1/2} \rho_t'\, dt \right\|^2 \\
            &\leq \int_0^1 \|(L_{w^{(\mu)}}^{\dagger})^{1/2} \rho_t'\|^2 \,dt\\
            &= \int_0^1 (\rho_t')^\top L_{w^{(\mu)}}^\dagger (\rho_t') \,dt,
        \end{align*}
    so that
        \begin{align*}
            \frac{1}{1 + C'|\varepsilon| + C''_{\eta}} \int_0^1 (\rho_t')^\top L_{w^{(\mu)}}^\dagger (\rho_t') dt &\geq \frac{1}{1 + C'|\varepsilon| + C''_{\eta}} \left(\int_0^1 \rho_t' dt\right)^\top L_{w^{(\mu)}}^{\dagger} \left(\int_0^1 \rho_t' dt\right)\\
            &= \frac{1}{1 + C'|\varepsilon| + C''_{\eta}} (\rho_1 - \rho_0)^\top L_{w^{(\mu)}}^{\dagger} (\rho_1 - \rho_0) \\
            &= \frac{\varepsilon^2}{1 + C'|\varepsilon| + C''_{\eta}} h^\top L_{w^{(\mu)}}^{\dagger} h.
        \end{align*}
    Since $\Maas(\mu,\mu+\varepsilon h)^2 \leq \mathcal{A}((\psi_t),(\rho_t)) \leq \Maas(\mu,\mu+\varepsilon h)^2 + \eta$, the previous inequality yields
        \begin{align*}
            \Maas(\mu, \mu+\varepsilon \, h)^2 + \eta &\geq \frac{\varepsilon^2}{1 + C'|\varepsilon| + C''_{\eta}} h^\top L_{w^{(\mu)}}^{\dagger} h.
        \end{align*}
    Letting $\eta\to 0$ gives the desired lower bound
        \begin{align*}
            \Maas(\mu, \mu+\varepsilon \, h)^2 &\geq \frac{\varepsilon^2}{1 + C'|\varepsilon|} h^\top L_{w^{(\mu)}}^{\dagger} h,
        \end{align*}
    and the theorem is proved.
\end{proof}

To conclude this subsection we briefly discuss the sharpness of the constants $C, C'$ appearing in the statement of~\Cref{thm:linearization-mu}. For the upper bound, our proof only uses the pointwise estimate
    \begin{align*}
        (1-|\varepsilon|C)\,\mu(x)\leq \rho_t(x)\leq (1+|\varepsilon|C)\,\mu(x),\qquad \rho_t:=\mu+t\varepsilon h,
    \end{align*}
which comes from the crude bound $|h(x)|\leq \|h\|_2\leq 1$. Thus, for a fixed direction $h$, one may replace
    \begin{align*}
        C=\left(\min_{x\in V}\mu(x)\right)^{-1}
    \end{align*}
by the sharper \emph{direction dependent} quantity
    \begin{align*}
        C_h:=\Bigl\|\frac{h}{\mu}\Bigr\|_{\ell^\infty(V)} =\max_{x\in V}\frac{|h(x)|}{\mu(x)}.
    \end{align*}
In particular, if $h$ is supported on vertices where $\mu$ is relatively large, then the upper bound improves accordingly. On the other hand, if one wants a statement uniform over all $h\in\mathbb{R}^V_0$ with $\|h\|_2\leq 1$, then some dependence on $(\min_x\mu(x))^{-1}$ is unavoidable, since positivity of $\mu+\varepsilon h$ already forces $|\varepsilon|$ to be at most of order $\min_x\mu(x)$ in the worst case.

Similarly, for the lower bound, our proof first inserts the coarse estimate
    \begin{align*}
        \Maas(\mu,\mu+\varepsilon h)^2 \leq 2\varepsilon^2\,h^\top L_{w^{(\mu)}}^\dagger h \leq \frac{2\varepsilon^2}{\lambda_2(L_{w^{(\mu)}})},
    \end{align*}
and then converts this into a uniform control on the deviation of an almost minimizing curve away from $\mu$. If one keeps the factor $h^\top L_{w^{(\mu)}}^\dagger h$ instead of bounding it by $\lambda_2(L_{w^{(\mu)}})^{-1}$, the same argument yields the sharper direction dependent constant
    \begin{align*}
        C_h' := \left(\min_{x\in V}\mu(x)\right)^{-1} \sqrt{2\,\Vol_w(G)\,h^\top L_{w^{(\mu)}}^\dagger h}.
    \end{align*}

Finally, the order of the error term can improve substantially for structured perturbations. The estimate in~\Cref{thm:linearization-mu} only gives a relative error of order $O(|\varepsilon|)$, because it tracks the first order variation of the metric tensor along an arbitrary affine path. In special cases this can be improved; this is the focus of the subsection which follows.

\subsection{Example: Bipartite graphs}\label{subsec:bipartite}

In general, the discrete transportation metric $\Maas(\cdot,\cdot)$ can be difficult to compute and describe in explicit terms. This makes a quantitative assessment of the quality of $\Resis(\cdot,\cdot)$ as a surrogate for $\Maas(\cdot,\cdot)$ difficult to assess in general. However, certain settings are more amenable to such an analysis, and the focus of this section is to describe one such example in detail.

Specifically, in the setting of a connected bipartite graph $G$, there exists a one parameter family of perturbations of the corresponding stationary distribution $\pi$ for which exact quantitative guarantees on the linearization of $\Maas(\cdot,\cdot)$ near $\pi$ can be provided. Moreover, the resulting approximation is stronger than the general bound provided by~\Cref{thm:linearization-mu}, in the sense that $\Resis(\cdot,\cdot)$ captures $\Maas(\cdot,\cdot)^2$ up to third order as the perturbation size tends to zero. 

Let $G=(V,E,w)$ denote a fixed connected bipartite graph with a corresponding bipartition $V=V_1\cup V_2$. We define the \emph{sign function} $\zeta\in\mathbb{R}^{V}$ according to
    \begin{align*}
        \zeta(x) &= \begin{cases}
            1&\text{ if }x\in V_1\\
            -1&\text{ if }x\in V_2
        \end{cases},\quad x\in V.
    \end{align*}
We consider the following single parameter family of perturbations of $\pi$:
    \begin{align}\label{eq:bipartite-perturbation}
        \pi^{(s)} &:= \pi\,\odot\,(1+s\zeta),\quad -1<s<1,
    \end{align}
where $\odot$ denotes componentwise multiplication of vectors. In plain terms, $\pi^{(s)}$ describes a perturbation of the stationary distribution in which the mass on $V_1$ is multiplied by a factor of $1+s$ and the mass on $V_2$ is multiplied by a factor of $1-s$. We illustrate this setup in~\Cref{fig:bipartite-perturbation}.

\begin{figure}[t]
    \centering

    \begin{subfigure}[t]{0.32\textwidth}
        \centering
        \begin{tikzpicture}[x=0.82cm,y=0.72cm]
            \coordinate (x1) at (0,5);
            \coordinate (x2) at (0,4);
            \coordinate (x3) at (0,3);
            \coordinate (x4) at (0,2);
            \coordinate (x5) at (0,1);
            \coordinate (x6) at (0,0);

            \coordinate (y1) at (3.6,5);
            \coordinate (y2) at (3.6,4);
            \coordinate (y3) at (3.6,3);
            \coordinate (y4) at (3.6,2);
            \coordinate (y5) at (3.6,1);
            \coordinate (y6) at (3.6,0);

            \draw[line width=1pt] (x1) -- (y1);
            \draw[line width=1pt] (x1) -- (y2);

            \draw[line width=1pt] (x2) -- (y2);
            \draw[line width=1pt] (x2) -- (y3);
            \draw[line width=1pt] (x2) -- (y4);

            \draw[line width=1pt] (x3) -- (y3);
            \draw[line width=1pt] (x3) -- (y4);
            \draw[line width=1pt] (x3) -- (y5);

            \draw[line width=1pt] (x4) -- (y4);
            \draw[line width=1pt] (x4) -- (y5);

            \draw[line width=1pt] (x5) -- (y5);

            \draw[line width=1pt] (x6) -- (y5);
            \draw[line width=1pt] (x6) -- (y6);

            \node[circle, fill=blue!50, minimum size=9pt, inner sep=0pt] at (x1) {};
            \node[circle, fill=blue!70, minimum size=9pt, inner sep=0pt] at (x2) {};
            \node[circle, fill=blue!70, minimum size=9pt, inner sep=0pt] at (x3) {};
            \node[circle, fill=blue!50, minimum size=9pt, inner sep=0pt] at (x4) {};
            \node[circle, fill=blue!30, minimum size=9pt, inner sep=0pt] at (x5) {};
            \node[circle, fill=blue!50, minimum size=9pt, inner sep=0pt] at (x6) {};

            \node[circle, fill=blue!30, minimum size=9pt, inner sep=0pt] at (y1) {};
            \node[circle, fill=blue!50, minimum size=9pt, inner sep=0pt] at (y2) {};
            \node[circle, fill=blue!50, minimum size=9pt, inner sep=0pt] at (y3) {};
            \node[circle, fill=blue!70, minimum size=9pt, inner sep=0pt] at (y4) {};
            \node[circle, fill=blue!90, minimum size=9pt, inner sep=0pt] at (y5) {};
            \node[circle, fill=blue!30, minimum size=9pt, inner sep=0pt] at (y6) {};

            \node at (0,5.8) {$V_1$};
            \node at (3.6,5.8) {$V_2$};
        \end{tikzpicture}
        \caption{Stationary measure $\pi$.}
    \end{subfigure}
    \hfill
    \begin{subfigure}[t]{0.32\textwidth}
        \centering
        \begin{tikzpicture}[x=0.82cm,y=0.72cm]
            \coordinate (x1) at (0,5);
            \coordinate (x2) at (0,4);
            \coordinate (x3) at (0,3);
            \coordinate (x4) at (0,2);
            \coordinate (x5) at (0,1);
            \coordinate (x6) at (0,0);

            \coordinate (y1) at (3.6,5);
            \coordinate (y2) at (3.6,4);
            \coordinate (y3) at (3.6,3);
            \coordinate (y4) at (3.6,2);
            \coordinate (y5) at (3.6,1);
            \coordinate (y6) at (3.6,0);

            \draw[line width=1pt] (x1) -- (y1);
            \draw[line width=1pt] (x1) -- (y2);

            \draw[line width=1pt] (x2) -- (y2);
            \draw[line width=1pt] (x2) -- (y3);
            \draw[line width=1pt] (x2) -- (y4);

            \draw[line width=1pt] (x3) -- (y3);
            \draw[line width=1pt] (x3) -- (y4);
            \draw[line width=1pt] (x3) -- (y5);

            \draw[line width=1pt] (x4) -- (y4);
            \draw[line width=1pt] (x4) -- (y5);

            \draw[line width=1pt] (x5) -- (y5);

            \draw[line width=1pt] (x6) -- (y5);
            \draw[line width=1pt] (x6) -- (y6);

            \node[circle, fill=blue!80, minimum size=9pt, inner sep=0pt] at (x1) {};
            \node[circle, fill=blue!95, minimum size=9pt, inner sep=0pt] at (x2) {};
            \node[circle, fill=blue!95, minimum size=9pt, inner sep=0pt] at (x3) {};
            \node[circle, fill=blue!80, minimum size=9pt, inner sep=0pt] at (x4) {};
            \node[circle, fill=blue!55, minimum size=9pt, inner sep=0pt] at (x5) {};
            \node[circle, fill=blue!80, minimum size=9pt, inner sep=0pt] at (x6) {};

            \node[circle, fill=blue!8, minimum size=9pt, inner sep=0pt] at (y1) {};
            \node[circle, fill=blue!12, minimum size=9pt, inner sep=0pt] at (y2) {};
            \node[circle, fill=blue!12, minimum size=9pt, inner sep=0pt] at (y3) {};
            \node[circle, fill=blue!18, minimum size=9pt, inner sep=0pt] at (y4) {};
            \node[circle, fill=blue!24, minimum size=9pt, inner sep=0pt] at (y5) {};
            \node[circle, fill=blue!8, minimum size=9pt, inner sep=0pt] at (y6) {};

            \node at (0,5.8) {$V_1$};
            \node at (3.6,5.8) {$V_2$};

        \end{tikzpicture}
        \caption{Perturbed measure $\pi^{(3/4)}$.}
    \end{subfigure}
    \hfill
    \begin{subfigure}[t]{0.32\textwidth}
        \centering
        \begin{tikzpicture}[x=0.82cm,y=0.72cm]
            \coordinate (x1) at (0,5);
            \coordinate (x2) at (0,4);
            \coordinate (x3) at (0,3);
            \coordinate (x4) at (0,2);
            \coordinate (x5) at (0,1);
            \coordinate (x6) at (0,0);

            \coordinate (y1) at (3.6,5);
            \coordinate (y2) at (3.6,4);
            \coordinate (y3) at (3.6,3);
            \coordinate (y4) at (3.6,2);
            \coordinate (y5) at (3.6,1);
            \coordinate (y6) at (3.6,0);

            \draw[line width=1pt] (x1) -- (y1);
            \draw[line width=1pt] (x1) -- (y2);

            \draw[line width=1pt] (x2) -- (y2);
            \draw[line width=1pt] (x2) -- (y3);
            \draw[line width=1pt] (x2) -- (y4);

            \draw[line width=1pt] (x3) -- (y3);
            \draw[line width=1pt] (x3) -- (y4);
            \draw[line width=1pt] (x3) -- (y5);

            \draw[line width=1pt] (x4) -- (y4);
            \draw[line width=1pt] (x4) -- (y5);

            \draw[line width=1pt] (x5) -- (y5);

            \draw[line width=1pt] (x6) -- (y5);
            \draw[line width=1pt] (x6) -- (y6);

            \node[circle, fill=blue!18, minimum size=9pt, inner sep=0pt] at (x1) {};
            \node[circle, fill=blue!24, minimum size=9pt, inner sep=0pt] at (x2) {};
            \node[circle, fill=blue!24, minimum size=9pt, inner sep=0pt] at (x3) {};
            \node[circle, fill=blue!18, minimum size=9pt, inner sep=0pt] at (x4) {};
            \node[circle, fill=blue!12, minimum size=9pt, inner sep=0pt] at (x5) {};
            \node[circle, fill=blue!18, minimum size=9pt, inner sep=0pt] at (x6) {};

            \node[circle, fill=blue!55, minimum size=9pt, inner sep=0pt] at (y1) {};
            \node[circle, fill=blue!80, minimum size=9pt, inner sep=0pt] at (y2) {};
            \node[circle, fill=blue!80, minimum size=9pt, inner sep=0pt] at (y3) {};
            \node[circle, fill=blue!95, minimum size=9pt, inner sep=0pt] at (y4) {};
            \node[circle, fill=blue!95, minimum size=9pt, inner sep=0pt] at (y5) {};
            \node[circle, fill=blue!55, minimum size=9pt, inner sep=0pt] at (y6) {};

            \node at (0,5.8) {$V_1$};
            \node at (3.6,5.8) {$V_2$};

        \end{tikzpicture}
        \caption{Perturbed measure $\pi^{(-3/4)}$.}
    \end{subfigure}

    \caption{A bipartite graph and the perturbation family
    $\pi^{(s)}=\pi\odot(1+s\zeta)$ from \cref{eq:bipartite-perturbation}. In (a), node color opacity at vertex $x\in V=V_1\cup V_2$ is proportional to the stationary distribution $\pi(x)$, equivalently to degree since $\pi(x)=\Deg_x/\Vol_w(G)$. In (b) and (c), the same graph is shown with color opacity proportional to $\pi^{(3/4)}(x)$ and $\pi^{(-3/4)}(x)$, respectively, illustrating how the perturbation shifts mass toward $V_1$ or $V_2$.}
    \label{fig:bipartite-perturbation}
\end{figure}

Recall that at a fixed reference measure $\mu\in\Prob_\star(V)$, the discrete transportation distance metric tensor is determined by the edge weights
    \begin{align*}
        w_{xy}^{(\mu)} = w_{xy}\,\theta\left(\frac{\mu(x)}{\pi(x)}, \frac{\mu(y)}{\pi(y)}\right),\quad \{x,y\}\in E.
    \end{align*}
In the case of~\cref{eq:bipartite-perturbation}, this reads
    \begin{align*}
        w_{xy}^{(\pi\odot(1+s\zeta))} = w_{xy}\,\theta(1+s, 1-s)
    \end{align*}
for each $\{x,y\}\in E$. In other words, the weights underwriting the corresponding metric tensor change only by a global factor $\theta(1+s, 1-s)$ as $s$ varies. 

The main result of this section is given below.

\begin{theorem}\label{thm:bipartite-third-order}
    Let $G=(V,E,w)$ be a connected bipartite graph and let $\pi$ denote the corresponding stationary distribution of the simple random walk on $V$. Let $\pi^{(s)}$ be as in~\cref{eq:bipartite-perturbation}. Then for $-1< a < b < 1$, as $(a,b)\to (0,0)$, it holds
        \begin{align*}
            \Maas(\pi^{(a)},\pi^{(b)})^2 &= \Resis(\pi^{(a)},\pi^{(b)}) + \frac{(b-a)^2(a^2+ab+b^2)}{18\,\Vol_w(G)} + O((|a|+|b|)^6).
        \end{align*}
    In particular, the resistance distance $\Resis(\pi^{(a)},\pi^{(b)})$ captures $\Maas(\pi^{(a)},\pi^{(b)})^2$ up to third order as $(a,b)\to(0,0)$.
\end{theorem}

Below we state and prove~\Cref{prop:bipartite-explicit} and~\Cref{prop:bipartite-explicit-resis}, which compute, respectively, $\Maas(\pi^{(a)},\pi^{(b)})^2$ and $\Resis(\pi^{(a)},\pi^{(b)})$ in exact terms. With these in hand,~\Cref{thm:bipartite-third-order} follows from Taylor's theorem and a short computation; we supply the details at the end of this section for completeness.

\begin{proposition}\label{prop:bipartite-explicit}
    Let $G=(V_1\cup V_2, E, w)$ be a connected bipartite graph. Let $\pi$ denote the corresponding stationary distribution of the simple random walk on $V$ and let $\pi^{(s)}$ be as in~\cref{eq:bipartite-perturbation}. For brevity, if $-1 < u <1$, write $\theta(u) := \theta(1+u, 1-u)$. Let $-1 < a < b < 1$ be fixed. Then
        \begin{align}\label{eq:bipartite-maas}
            \Maas(\pi^{(a)},\pi^{(b)}) &= \frac{1}{\sqrt{2\Vol_w(G)}} \left| \int_a^{b} \frac{du}{\sqrt{\theta(u)}}\,\right|.
        \end{align}
\end{proposition}

\begin{proof}
    First we prove the upper bound in~\cref{eq:bipartite-maas}:
        \begin{align}\label{eq:bipartite-maas-upper}
            \Maas(\pi^{(a)},\pi^{(b)}) &\leq \frac{1}{\sqrt{2\Vol_w(G)}} \left| \int_a^{b} \frac{du}{\sqrt{\theta(u)}}\,\right|.
        \end{align}
    Let $\varphi\in\mathcal{C}^{1}([0,1];(-1,1))$ satisfy $\varphi(0) = a$ and $\varphi(1) = b$. Define, for $0\le t\le 1$,
        \begin{align*}
            \rho_t &:= \pi^{(\varphi(t))},\quad \psi_t := \frac{\varphi'(t)}{2\Vol_w(G) \theta(\varphi(t))}\zeta.
        \end{align*}
    We claim that $(\rho_t,\psi_t)$ is feasible for $\Maas(\pi^{(a)},\pi^{(b)})$ as in~\cref{eq:defn-maas-distance}. Clearly $(\rho_t)\in\mathcal{C}^{1}([0,1];\Prob_\star(V))$, $\rho_0 = \pi^{(a)}$, $\rho_1 = \pi^{(b)}$, and $(\psi_t)\in \mathcal{B}([0,1];\mathbb{R}^{V})$ by construction. For the continuity equation, first observe that
        \begin{align*}
            \frac{d}{dt} \rho_t &= \frac{d}{dt} \pi^{(\varphi(t))} = \varphi'(t)  \pi\,\odot\,\zeta.
        \end{align*}
    In the notation of~\Cref{thm:linearization-mu}, for $-1<s<1$, the weights $w_{xy}^{\pi^{(s)}} =: w_{xy}^{(s)}$ can be computed as
        \begin{align*}
            w_{xy}^{(s)} &= w_{xy} \, \theta\left(\frac{\pi^{(s)}(x)}{\pi(x)}, \frac{\pi^{(s)}(y)}{\pi(y)}\right) = w_{xy}\, \theta(1+s, 1-s),
        \end{align*}
    for each $\{x, y\}\in E$. For brevity, write $\theta(s):= \theta(1+s, 1-s)$. Then it follows that
        \begin{align}\label{eq:bipartite-Lpsi}
            L_{w^{(\varphi(t))}}\psi_t &= \frac{\varphi'(t)}{2\Vol_w(G) \theta(\varphi(t))} L_{w^{(\varphi(t))}}\zeta.
        \end{align}
    We compute, for $x\in V_1$,
        \begin{align*}
            L_{w^{(\varphi(t))}}\zeta(x) &= \theta(\varphi(t)) L_w \zeta(x) = \theta(\varphi(t)) \sum_{y\sim x} w_{xy}(\zeta(x)-\zeta(y)) = 2 \theta(\varphi(t)) \Deg_x \zeta(x),
        \end{align*}
    and similarly for $x\in V_2$. Thus~\cref{eq:bipartite-Lpsi} becomes
        \begin{align*}
            L_{w^{(\varphi(t))}}\psi_t & = \varphi'(t)\, \pi\odot \zeta = \frac{d}{dt}\rho_t,
        \end{align*}
    as desired. We compute the action functional as follows:
        \begin{align*}
            \mathcal{A}((\psi_t), (\rho_t)) &= \int_0^1 \psi_t^\top L_{w^{(\varphi(t))}} \psi_t\, dt = \frac{1}{2\Vol_w(G)}\int_0^1 \frac{\varphi'(t)^2}{\theta(\varphi(t))} \zeta^\top (\pi\odot \zeta) \, dt
        \end{align*}
    We observe that $ \zeta^\top (\pi\odot \zeta) = 1$, so 
        \begin{align*}
            \mathcal{A}((\psi_t), (\rho_t)) &= \frac{1}{2\Vol_w(G)}\int_0^1 \frac{\varphi'(t)^2}{\theta(\varphi(t))} \, dt.
        \end{align*}
    It remains to choose $\varphi$ judiciously. To this end, write
        \begin{align*}
            \Phi(s) &:= \int_0^{s} \frac{du}{\sqrt{\theta(u)}},\quad -1<s<1.
        \end{align*}    
    Then $\Phi(s)$ is strictly increasing in $s$ and admits an inverse on its image. Put
        \begin{align*}
            \varphi(t) &:= \Phi^{-1}((1-t)\Phi(a) + t\Phi(b)),\quad 0\le t\le 1.
        \end{align*}
    Then 
        \begin{align*}
            \dfrac{d}{dt} \Phi(\varphi(t)) &= \Phi(b) - \Phi(a),
        \end{align*}
    and separately,
        \begin{align*}
            \dfrac{d}{dt} \Phi(\varphi(t)) &= \frac{\varphi'(t)}{\sqrt{\theta(\varphi(t))}},
        \end{align*}
    so that $\frac{\varphi'(t)}{\sqrt{\theta(\varphi(t))}} = \Phi(b) - \Phi(a)$ for all $0\le t\le 1$. Thus
        \begin{align*}
            \mathcal{A}((\psi_t), (\rho_t)) &= \frac{1}{2\Vol_w(G)}\, \int_0^1 \frac{\varphi'(t)^2}{\theta(\varphi(t))} \, dt = \frac{1}{2\Vol_w(G)}\, (\Phi(b) - \Phi(a))^2 = \frac{1}{2\Vol_w(G)}\, \left| \int_a^{b} \frac{du}{\sqrt{\theta(u)}}\right|^2.
        \end{align*}
    Thus~\cref{eq:bipartite-maas-upper} is proved. Next we prove the lower bound
        \begin{align}\label{eq:bipartite-maas-lower}
            \Maas(\pi^{(a)},\pi^{(b)}) &\geq \frac{1}{\sqrt{2\Vol_w(G)}} \left| \int_a^{b} \frac{du}{\sqrt{\theta(u)}}\,\right|.
        \end{align}
    To this end, let $((\rho_t),(\psi_t))$ be a feasible pair for $\Maas(\pi^{(a)},\pi^{(b)})$ in the sense of ~\cref{eq:defn-maas-distance}. Define the function $\sigma(t) \in \mathcal{C}^{1}([0,1])$ via
        \begin{align*}
            \sigma(t) &:= \sum_{x\in V} \zeta(x) \rho_t(x)\,\quad 0\le t\le 1.
        \end{align*}
    Then we have that, by feasibility of $\psi_t$, it holds
        \begin{align*}
            \frac{d}{dt}\sigma(t) &= \sum_{x\in V} \zeta(x) \frac{d}{dt}\rho_t(x) =  \sum_{x\in V} \zeta(x) L_{\widehat{w}^{(t)}} \psi_t(x),
        \end{align*}
    where we define the edge weights
        \begin{align*}
            \widehat{w}^{(t)}_{xy} &:= w_{xy}\theta\left(\frac{\rho_t(x)}{\pi(x)}, \frac{\rho_t(y)}{\pi(y)}\right).
        \end{align*}
    We compute
        \begin{align*}
            \sum_{x\in V} \zeta(x) L_{\widehat{w}^{(t)}} \psi_t(x) &= \sum_{x\in V} \zeta(x) \sum_{y\sim x} \widehat{w}_{xy}^{(t)} (\psi_t(x) - \psi_t(y)) \\
            &= \sum_{x\in V_1} \sum_{y\sim x} \widehat{w}_{xy}^{(t)} (\psi_t(x) - \psi_t(y)) - \sum_{x\in V_2} \sum_{y\sim x} \widehat{w}_{xy}^{(t)} (\psi_t(x) - \psi_t(y))\\
            &= \sum_{(x,y)\in E'} 2 \widehat{w}_{xy}^{(t)} (\psi_t(x) - \psi_t(y)),
        \end{align*}
    where $E'$ is the orientation of $E$ consisting of edges directed from $V_1$ to $V_2$. By the Cauchy-Schwarz inequality, we have
        \begin{align}\label{eq:bipartite-cauchy-schwarz}
            \left|\frac{d}{dt}\sigma(t)\right|^2 &\leq 4 \left(\sum_{(x,y)\in E'} \widehat{w}_{xy}^{(t)} (\psi_t(x) - \psi_t(y))^2\right) \left(\sum_{(x,y)\in E'} \widehat{w}_{xy}^{(t)}\right)
        \end{align}
    We estimate the rightmost factor as follows. First recall that $(s,t)\mapsto \theta(s, t)$ is jointly concave in $s,t$ on $(0,\infty)\times(0,\infty)$. Therefore by Jensen's inequality,
        \begin{align*}
            \sum_{(x,y)\in E'} \frac{w_{xy}}{\frac{1}{2}\Vol_w(G)} \theta\left(\frac{\rho_t(x)}{\pi(x)}, \frac{\rho_t(y)}{\pi(y)}\right) &\leq \theta\left( \frac{ \sum_{(x,y)\in E'} w_{xy} \frac{\rho_t(x)}{\pi(x)}}{\frac{1}{2}\Vol_w(G)},  \frac{ \sum_{(x,y)\in E'} w_{xy} \frac{\rho_t(y)}{\pi(y)}}{\frac{1}{2}\Vol_w(G)} \right)\\
            &= \theta\left(  2\sum_{(x,y)\in E'} \frac{w_{xy}}{\Deg_x}\rho_t(x), 2\sum_{(x,y)\in E'} \frac{w_{xy}}{\Deg_y}\rho_t(y) \right)\\
            &= \theta\left( 2\sum_{x\in V_1} \rho_t(x), 2\sum_{y\in V_2} \rho_t(y) \right) = \theta(\sigma(t)).
        \end{align*}
    Therefore,
        \begin{align*}
            \sum_{(x,y)\in E'} \widehat{w}_{xy}^{(t)} &\leq \frac{\Vol_w(G)}{2}\,\theta(\sigma(t)).
        \end{align*}
    Combining the preceding with~\cref{eq:bipartite-cauchy-schwarz} and integrating, it follows that
        \begin{align*}
            \mathcal{A}((\psi_t),(\rho_t)) &\geq \frac{1}{2\Vol_w(G)} \int_0^1 \frac{|\sigma'(t)|^2}{\theta(\sigma(t))} dt.
        \end{align*}
    Since $\Phi'(s) = \theta(s)^{-1/2}$, one has
        \begin{align*}
            \frac{d}{dt}\Phi(\sigma(t)) = \frac{\sigma'(t)}{\sqrt{\theta(\sigma(t))}}. 
        \end{align*}
    By the Cauchy-Schwarz inequality, we have
        \begin{align*}
            \int_0^1 \frac{|\sigma'(t)|^2}{\theta(\sigma(t))} dt &= \int_0^1 \left|\frac{d}{dt}\Phi(\sigma(t))\right|^2 dt \geq \left|\int_0^1 \frac{d}{dt}\Phi(\sigma(t)) dt\right|^2 = |\Phi(\sigma(1)) - \Phi(\sigma(0))|^2,
        \end{align*}
    since $\Phi$ is strictly increasing. Finally, since $\sigma(0) = \sum_{x\in V} \zeta(x)\pi^{(a)}(x) = a$ and $\sigma(1) = \sum_{x\in V} \zeta(x)\pi^{(b)}(x) = b$, we have
        \begin{align*}
            |\Phi(\sigma(1)) - \Phi(\sigma(0))|^2 = |\Phi(b) - \Phi(a)|^2 = \left|\int_a^{b} \frac{du}{\sqrt{\theta(u)}}\right|^2.
        \end{align*}
    Taking the infimum over all feasible pairs $((\rho_t),(\psi_t))$ yields the lower bound~\cref{eq:bipartite-maas-lower}, and the proof is complete.
\end{proof}

\begin{proposition}\label{prop:bipartite-explicit-resis}
    Let $G=(V_1\cup V_2, E, w)$ be a connected bipartite graph. Let $\pi$ denote the corresponding stationary distribution of the simple random walk on $V$ and let $\pi^{(s)}$ be as in~\cref{eq:bipartite-perturbation}. Let $-1 < a < b < 1$ be fixed. Then
        \begin{align}\label{eq:bipartite-maas-resistance}
            \Resis(\pi^{(a)},\pi^{(b)}) &= \frac{(a-b)^2}{ 2\Vol_w(G)}.
        \end{align}
\end{proposition}

\begin{proof}
    Write
        \begin{align*}
            h := \pi \odot \zeta \in \mathbb{R}^{V}.
        \end{align*}
    Since $\pi(V_1)=\pi(V_2)=\frac12$, it follows that
        \begin{align*}
            \sum_{x\in V} h(x) = \sum_{x\in V} \pi(x)\zeta(x) = \pi(V_1)-\pi(V_2)=0,
        \end{align*}
    and hence $h\in\mathbb{R}^{V}_0$. By definition of $\pi^{(s)}$,
        \begin{align}\label{eq:bipartite-resis-difference}
            \pi^{(a)} - \pi^{(b)} = (a-b)\,h.
        \end{align}
    As in the proof of~\Cref{prop:bipartite-explicit}, it holds
        \begin{align}\label{eq:bipartite-resis-Lzeta}
            (L_{w}\zeta)(x) &= 2\Vol_w(G)\,h(x).
        \end{align}
    Using~\Cref{defn:resistance-distance-measures,eq:bipartite-resis-difference}, we obtain
        \begin{align*}
            \Resis(\pi^{(a)},\pi^{(b)})
            &= (\pi^{(a)}-\pi^{(b)})^\top L_w^\dagger (\pi^{(a)}-\pi^{(b)}) \\
            &= (a-b)^2 h^\top L_w^\dagger h.
        \end{align*}
    It therefore remains to evaluate $h^\top L_w^\dagger h$. By~\cref{eq:bipartite-resis-Lzeta},
        \begin{align*}
            h^\top L_w^\dagger h
            &= \frac{1}{2\Vol_w(G)}\, h^\top L_w^\dagger L_w \zeta.
        \end{align*}
    Since $h\in\mathbb{R}^{V}_0$ and $L_w L_w^\dagger$ is the orthogonal projection onto $\mathbb{R}^{V}_0$, it holds
        \begin{align*}
            \frac{1}{2\Vol_w(G)}\, h^\top L_w^\dagger L_w \zeta &= \frac{1}{2\Vol_w(G)}\, h^\top \zeta =  \frac{1}{2\Vol_w(G)} (\zeta\,\odot\,\pi)^\top\zeta = \frac{1}{2\Vol_w(G)}.
        \end{align*}
    Therefore
        \begin{align*}
            \Resis(\pi^{(a)},\pi^{(b)}) &= \frac{(a-b)^2}{2\Vol_w(G)},
        \end{align*}
    as claimed.
\end{proof}

\begin{proof}[Proof of~\Cref{thm:bipartite-third-order}]
    Once again we write $\theta(u) := \theta(1-u, 1+u)$ for $-1<u<1$. One computes the Taylor expansion
        \begin{align*}
            \theta(u)=\theta(1+u,1-u)=\frac{u}{\operatorname{arctanh}(u)}
            = 1-\frac{u^2}{3}-\frac{4u^4}{45}+O(u^6),
        \end{align*}
    which is valid as $u\to 0$ (in particular the Taylor series has an interval of convergence $-1<u<1$). By the binomial power series we then have
        \begin{align*}
            \frac{1}{\sqrt{\theta(u)}}
            &= 1+\frac{u^2}{6}+\frac{31u^4}{360}+O(u^6).
        \end{align*}
    Thus, by~\Cref{prop:bipartite-explicit} and by~\Cref{prop:bipartite-explicit-resis}, it holds
        \begin{align*}
            \Maas(\pi^{(a)},\pi^{(b)})^2 &= \frac{1}{2\Vol_w(G)} \left| \int_a^{b} \frac{du}{\sqrt{\theta(u)}}\,\right|^2\\
            &= \frac{1}{2\Vol_w(G)} \left| \int_a^{b} \left(1+\frac{u^2}{6}+\frac{31u^4}{360}+O(u^6)\right) du\,\right|^2\\
            &= \frac{1}{2\Vol_w(G)}\left((b-a)+\frac{b^3-a^3}{18}+O((|a|+|b|)^5)\right)^2\\
            &= \frac{(b-a)^2}{2\,\Vol_w(G)} + \frac{(b-a)^2(a^2+ab+b^2)}{18\,\Vol_w(G)} + O((|a|+|b|)^6)\\
            &= \Resis(\pi^{(a)},\pi^{(b)}) + \frac{(b-a)^2(a^2+ab+b^2)}{18\,\Vol_w(G)} + O((|a|+|b|)^6),
        \end{align*}
    again as $(a,b) \to 0$.
\end{proof}

\section{Characterizations of resistance distance}\label{sec:characterizations}

In~\Cref{sec:linearization}, we showed that near a reference measure, the discrete transportation distance admits a second-order expansion governed by a Laplacian pseudoinverse, and that in the stationary case $\mu=\pi$ the corresponding inner product is exactly $\langle f,g\rangle_{L_w^\dagger}=f^\top L_w^\dagger g$ on $\mathbb{R}^V_0$. In light of this, it is natural to regard the interior of the probability simplex as a flat Riemannian manifold whose metric tensor is induced by $L_w^\dagger$, and, in doing so, to explore the dynamical and variational statements that become available in the resulting geometry. The goal of this section is to formalize this viewpoint, and then develop a number of equivalent characterizations of the resulting geodesic distance.

We begin by restating the definition of the so-called resistance manifold.

\begin{definition}[Resistance manifold]\label{defn:resistance-manifold}
    Let $G=(V,E,w)$ be a connected graph, and for $\mu,\nu\in\Prob_\star(V)$, write
        \begin{align}\label{eq:distance-resis}
            d_\Resis(\mu,\nu) &:= \sqrt{\Resis(\mu,\nu)} = \sqrt{(\mu-\nu)^\top L_w^\dagger (\mu-\nu)}.
        \end{align}
    The \emph{resistance manifold} is the metric space $(\Prob_\star(V), d_\Resis)$. 
\end{definition}

We remark that $(\Prob_\star(V), d_\Resis)$ is naturally a Riemannian manifold: $\Prob_\star(V)$ is an open subset of the affine hyperplane $\mathcal{H}=\{\rho\in\mathbb{R}^V:\sum_{x\in V}\rho(x)=1\}$ which can be equipped with the Riemannian metric tensor $\mathsf{g}^\Resis$ defined by
    \begin{align}\label{eq:resistance-inner-product}
        \mathsf{g}^\Resis_\rho(f, g) := \langle f, g\rangle_{L_{w}^\dagger},\qquad \rho\in\Prob_\star(V),\ \ f, g\in T_\rho\Prob_\star(V)\cong \mathbb{R}^V_0,
    \end{align}
and the Riemannian distance induced by $\mathsf{g}^\Resis$ coincides with $d_\Resis$. Moreover, since $\mathsf{g}^\Resis_\rho$ does not depend on $\rho$, the manifold is flat, and its geodesics are precisely the Euclidean line segments in $\mathcal{H}$. We state this as a lemma below without proof.

\begin{lemma}\label{lem:resistance-geodesics}
    Let $\rho_0,\rho_1\in\Prob_\star(V)$, and set $\rho_t := (1-t)\rho_0+t\rho_1$. Then $(\rho_t)_{t\in[0,1]}\subseteq\Prob_\star(V)$ is the unique constant-speed geodesic joining $\rho_0$ and $\rho_1$ in the resistance manifold.
\end{lemma}

Thus, in the remaining two subsections, we aim to provide several characterizations the metric $\Resis(\cdot,\cdot)$ and in the process prove~\Cref{thm:characterizations}. Many of these are based on classical results for resistance distance between vertices, and their proofs amount to generalizations in this regard. In doing so we also highlight some related facts which may be of interest to the reader.~\Cref{subsec:variational-characterizations} focuses on variational characterizations of $\Resis(\cdot,\cdot)$, while~\Cref{subsec:combinatorial-characterizations} focuses on combinatorial characterization of $\Resis(\cdot,\cdot)$.

In the interest of consistency with the results in~\Cref{sec:linearization} as well as~\Cref{sec:gradient-flows} to follow, in the definition of the resistance manifold we consider only measures in $\Prob_\star(V)$. However, the characterizations of $\Resis(\cdot,\cdot)$ which follow remain valid if one or more entries of the measure(s) under consideration vanish, so these are stated on $\Prob(V)$ to achieve fuller generality.

\subsection{Variational characterizations of resistance distance}\label{subsec:variational-characterizations}

In this subsection we prove the first half of~\Cref{thm:characterizations}, which consists of various variational descriptions of resistance distance between measures.

First we recall some background concerning so-called \emph{Beckmann problems}, beginning with the $p$-Wasserstein metric.

\begin{definition}[$p$-Wasserstein metric]\label{defn:wasserstein}
    Let $G=(V, E, w)$ be a connected graph and let $\mu,\nu\in\Prob(V)$ and $p\geq 1$ be fixed. We define their \emph{$p$-Wasserstein distance}, denoted $\Wass_{p}(\mu,\nu)$, according to the linear program
        \begin{align*}
            \Wass_{p}(\mu,\nu)^p &= \inf_{\gamma\in\Gamma(\mu,\nu)} \sum_{x,y\in V} \Dist_{w}(x,y)^p \, \gamma(x,y),
        \end{align*}
    where $\Gamma(\mu,\nu)$ is the set of couplings between $\mu$ and $\nu$, i.e.,
        \begin{align*}
            \Gamma(\mu,\nu) &= \left\{\gamma\in\mathbb{R}^{V\times V}:\, \gamma(x,y)\geq 0,\, \sum_{y\in V} \gamma(x,y) = \mu(x),\, \sum_{x\in V} \gamma(x,y) = \nu(y),\, x,y\in V \right\},
        \end{align*}
    and where, for $x, y\in V$, $\Dist_{w}(x, y)$ is the weighted shortest path distance between $x$ and $y$.
\end{definition}

In this subsection, the graph metric in~\Cref{defn:wasserstein} is understood as the shortest path metric with edge lengths given by the values $w_{xy}$.

In the setting $p=1$, the Wasserstein distance admits several equivalent characterizations, including the so-called \emph{Beckmann formulation}, which takes the following form:
    \begin{align}\label{eq:beckmann-w1}
        \Wass_{1}(\mu,\nu) &= \inf_{J\in\mathbb{R}^{E'}} \left\{\sum_{(x, y)\in E'} w_{xy} |J(x,y)|:\, B J = \mu - \nu \right\}.
    \end{align}
Here, $E'$ is an orientation of $E$, and $B\in\mathbb{R}^{V\times E'}$ is the associated \emph{oriented incidence matrix} with entries given by
    \begin{align*}
        B_{x, (u,v)} &= \begin{cases}
            1,& x=u,\\
            -1,& x=v,\\
            0,& \text{otherwise}.
        \end{cases}
    \end{align*}
The equivalence exhibited in~\cref{eq:beckmann-w1} can be proved by two applications of strong duality for linear programs; we refer the reader to~\cite[Ch. 6]{peyre2019computational} for a proof. In the case of resistance distance $\Resis(\mu,\nu)$,~\cref{eq:beckmann-w1} admits a natural extension.

\begin{theorem}\label{thm:beckmann}
    Let $G=(V,E,w)$ be a connected graph and let $\mu,\nu\in\Prob(V)$. Let $E'$ denote an orientation of $E$, and let $B\in\mathbb{R}^{V\times E'}$ denote the associated oriented incidence matrix. Then the following holds:
        \begin{align}\label{eq:2-beckmann}
            \Resis(\mu,\nu) &= \inf_{J\in\mathbb{R}^{E'}} \left\{\sum_{e\in E'} \frac{|J(e)|^2}{w_{e}}:\, B J = \mu - \nu \right\}.
        \end{align}
\end{theorem}

\begin{proof}
    Let $W\in\mathbb{R}^{E'\times E'}$ denote the diagonal matrix with entries $W_{e,e} = w_e$ for $e\in E'$. Then~\cref{eq:2-beckmann} amounts to a linearly constrained least squares problem: minimize $\|W^{-1/2} J\|_{2}^{2}$ subject to $ B J = \mu - \nu$. By changing variables we can recast this as minimizing $\|Z\|_{2}^{2}$ subject to $ B W^{1/2} Z = \mu - \nu$. The Moore--Penrose inverse of the matrix $(BW^{1/2})\in\mathbb{R}^{V\times E'}$ recovers this solution; therefore
        \begin{align*}
            \inf_{J\in\mathbb{R}^{E'}} \left\{\sum_{e\in E'} \frac{|J(e)|^2}{w_{e}}:\, B J = \mu - \nu \right\} &= \| (B W^{1/2})^\dagger (\mu-\nu) \|_2^2\\
            &=  (\mu-\nu)^\top ((B W^{1/2})^\dagger)^\top  (B W^{1/2})^\dagger (\mu-\nu)\\
            &=  (\mu-\nu)^\top ((B W^{1/2}) (B W^{1/2})^\top )^\dagger (\mu-\nu)\\
            &=  (\mu-\nu)^\top L_w^\dagger (\mu-\nu) = \Resis(\mu,\nu).
        \end{align*}
\end{proof}

Our second variational characterization of resistance distance between measures is a version of the well known Benamou--Brenier formula for $2$-Wasserstein distance (see~\cite{benamou2000computational}). 

\begin{theorem}\label{thm:benamou-brenier-graphs}
    Let $G=(V,E,w)$ be a connected graph and let $\mu,\nu\in\Prob(V)$. Let $E'$ denote an orientation of $E$, and let $B\in\mathbb{R}^{V\times E'}$ denote the associated oriented incidence matrix. Then it holds
        \begin{align}\label{eq:benamou-brenier-graph}
            \Resis (\mu,\nu) &= \inf \Bigl\{\int_0^1 \sum_{e\in E'} \frac{|J_t(e)|^2}{w_e} \, dt \,:\, \rho_t\in \mathcal{C}^1([0,1]; \mathbb{R}^{V}), J_t\in \mathcal{B}([0,1]; \mathbb{R}^{E'}),\\
            &\qquad\qquad \frac{d}{dt} \rho_t + B J_t = 0,\quad \rho_0 = \mu,\quad \rho_1 = \nu\Bigr\}.
        \end{align}
\end{theorem}

\begin{proof}
    We denote by $\mathcal{I}$ the infimum on the right-hand side of~\cref{eq:benamou-brenier-graph}. For the upper bound $\mathcal{I} \le \Resis (\mu,\nu)$, let $J^\ast \in \mathbb{R}^{E'}$ be the unique minimizer of the Beckmann problem as in~\Cref{thm:beckmann}:
        \begin{align*}
            J^\ast &= \arg\min \left\{ \sum_{e\in E'} \frac{|J(e)|^2}{w_e} \,:\, B J = \mu - \nu \right\}.
        \end{align*}
    We construct a feasible path by taking the linear interpolation $\rho_t = (1-t) \mu + t \nu $, $0\le t\le 1$. Then $\rho_t \in \mathcal{C}^1([0,1]; \mathbb{R}^{V})$, $\rho_0 = \mu$, and $\rho_1 = \nu$. For this path, we have $\frac{d}{dt}\rho_t = \nu - \mu$. Consequently, the constraint $\frac{d}{dt}\rho_t + BJ_t = 0$ becomes
        \begin{align*}
            \nu - \mu + B J_t &= 0,
        \end{align*}
    which gives $B J_t = \mu - \nu$ for all $t \in [0, 1]$. To satisfy the continuity equation, we take the constant curve $J_t = J^*$ (independent of $t$). With this choice, the objective functional evaluates to
        \begin{align*}
            \int_0^1 \sum_{e\in E'} \frac{|J_t(e)|^2}{w_e} \, dt &= \int_0^1 \sum_{e\in E'} \frac{|J^\ast(e)|^2}{w_e} \, dt = \sum_{e\in E'} \frac{|J^\ast(e)|^2}{w_e} = \Resis(\mu,\nu).
        \end{align*}
    This shows that $\mathcal{I} \le \Resis(\mu,\nu)$. For the lower bound $\mathcal{I} \ge \Resis(\mu,\nu)$, let $(\rho_t, J_t)$ be a given admissible pair. We have
        \begin{align*}
            \frac{d}{dt} \rho_t + B J_t &= 0, \quad \rho_0 = \mu, \quad \rho_1 = \nu.
        \end{align*}
    Integrating the continuity equation from $t=0$ to $t=1$, we obtain
        \begin{align*}
            \nu - \mu + \int_0^1 B J_t \, dt &= 0,
        \end{align*}
    so
        \begin{align*}
            B \left( \int_0^1 J_t \, dt \right) &= \mu - \nu.
        \end{align*}
    Let $\overline{J} := \int_0^1 J_t \, dt$. Then $\overline{J} \in \mathbb{R}^{E'}$ satisfies the constraint $B\overline{J} = \mu - \nu$ that appears in the Beckmann problem. By~\Cref{thm:beckmann} and the Cauchy--Schwarz inequality,
        \begin{align*}
            \Resis(\mu,\nu) &\le \sum_{e\in E'} \frac{|\overline{J}(e)|^2}{w_e} = \sum_{e\in E'} \frac{\left| \int_0^1 J_t(e) \, dt \right|^2}{w_e} \le \sum_{e\in E'} \frac{ \int_0^1 |J_t(e)|^2 \, dt}{w_e} = \int_0^1 \sum_{e\in E'} \frac{|J_t(e)|^2}{w_e} \, dt.
        \end{align*}
    Since this holds for every feasible pair $(\rho_t, J_t)$, taking the infimum over all such pairs yields
        \begin{align*}
            \Resis(\mu,\nu) &\le \mathcal{I}.
        \end{align*}
    The claim is proved.
\end{proof}

The variational formula~\cref{eq:benamou-brenier-graph} can also be viewed as a relaxation of that in the definition of the discrete transportation distance (cf.~\cref{eq:defn-maas-distance}), wherein the weights $\theta\left(\frac{\rho_t(x)}{\pi(x)}, \frac{\rho_t(y)}{\pi(y)}\right)$ appearing in the action functional are replaced by constant unit weights.

The third and final variational characterization of resistance distance between measures is obtained by recasting $\Resis(\mu,\nu)$ as a negative homogeneous Sobolev norm distance interpreted in the graph setting. This is accomplished as follows. First, on $\mathbb{R}^{V}_0$, we define the following inner product:
    \begin{align*}
        \langle f, g \rangle_{\dot{H}^{1}(V,w)} &= \sum_{\{x,y\}\in E} w_{xy} (f(x) - f(y))(g(x) - g(y)),\quad f,g\in\mathbb{R}^{V}_0.
    \end{align*}
We remark that $\langle f, g \rangle_{\dot{H}^{1}(V,w)}$ is obviously bilinear and symmetric, and from
    \begin{align}\label{eq:hot-dot-one-on-graphs}
        \langle f, f \rangle_{\dot{H}^{1}(V,w)} &= \sum_{\{x,y\}\in E} w_{xy} (f(x) - f(y))^2 = f^\top L_w f,
    \end{align}
it follows that $\langle \cdot, \cdot \rangle_{\dot{H}^{1}(V,w)}$ is positive definite on $\mathbb{R}^{V}_0$. Thus, $\langle \cdot, \cdot \rangle_{\dot{H}^{1}(V,w)}$ defines an inner product on $\mathbb{R}^{V}_0$. The notation $\dot{H}^{1}(V,w)$ is taken from~\cref{eq:hot-dot-one-on-graphs}, by viewing $f^\top L_w f = \|W^{1/2}B^\top f\|_2^2$ as a (weighted) norm of the gradient of $f\in\mathbb{R}^{V}_0$.

By passing over to the formal dual of the norm induced by $\langle \cdot, \cdot \rangle_{\dot{H}^{1}(V,w)}$, we obtain a negative homogeneous Sobolev-type norm $\|\cdot\|_{\dot{H}^{-1}(V,w)}$ as adapted to the graph setting. The following theorem asserts that the norm distance between measures in this space recovers resistance distance.

\begin{theorem}\label{thm:sobolev-characterization}
    Let $G=(V,E,w)$ be a connected graph and let $\mu,\nu\in\Prob(V)$. Then it holds that
        \begin{align*}
            \Resis(\mu,\nu) &= \|\mu - \nu\|_{\dot{H}^{-1}(V,w)}^2 := \sup_{f\in\mathbb{R}^{V}_0\setminus\{0\}} \frac{\langle f, \mu - \nu \rangle^2}{\langle f, f \rangle_{\dot{H}^{1}(V,w)}}.
        \end{align*}
\end{theorem}

\begin{proof}
    By the definition of the $\dot{H}^{-1}(V,w)$ norm and~\cref{eq:hot-dot-one-on-graphs}, we have
        \begin{align*}
            \|\mu - \nu\|_{\dot{H}^{-1}(V,w)}^2 &= \sup_{f\in\mathbb{R}^{V}_0\setminus\{0\}} \frac{(f^\top (\mu - \nu))^2}{f^\top L_w f} = \sup_{f\in\mathbb{R}^{V}_0\setminus\{0\}} \frac{(f^\top (\mu - \nu))^2}{(L_w^{1/2} f)^\top (L_w^{1/2} f)}.
        \end{align*}
    On $\mathbb{R}^{V}_0$, the matrix $L_w$ and its square root are positive definite. Therefore we may apply a change of variables, writing $g=L_w^{1/2} f$, so that $f = L_w^{\dagger/2} g$ on $\mathbb{R}^{V}_0$. Thus,
        \begin{align*}
            \sup_{f\in\mathbb{R}^{V}_0\setminus\{0\}} \frac{(f^\top (\mu - \nu))^2}{(L_w^{1/2} f)^\top (L_w^{1/2} f)} &= \sup_{g\in\mathbb{R}^{V}_0\setminus\{0\}} \frac{\left((L_w^{\dagger/2} g)^\top (\mu - \nu)\right)^2}{g^\top g} = \sup_{g\in\mathbb{R}^{V}_0\setminus\{0\}} \frac{\left(g^\top L_w^{\dagger/2} (\mu - \nu)\right)^2}{g^\top g}.
        \end{align*}
    By the Cauchy--Schwarz inequality, the last expression equals
        \begin{align*}
            \|L_w^{\dagger/2} (\mu - \nu)\|_2^2 &= (\mu - \nu)^\top L_w^{\dagger} (\mu - \nu) = \Resis(\mu,\nu),
        \end{align*}
    and the theorem is proved.
\end{proof}

\subsection{Combinatorial characterizations of resistance distance}\label{subsec:combinatorial-characterizations}

In this subsection we prove the second half of~\Cref{thm:characterizations}, which consists of two combinatorial formulations of resistance distance between measures.

The first such formulation concerns expected vertex hitting times. To this end, let $(X_t)_{t\geq0}$ denote the simple random walk on $V$. For $y\in V$, we define the \emph{hitting time} to reach $y$, denoted $T_y$, as the random variable
    \begin{align*}
        T_y &= \inf \left\{t\ge 0 \,:\, X_{t} = y\right\}.
    \end{align*}
For $x, y\in V$, the \emph{expected hitting time} from $x$ to $y$, denoted $\Hitting(x, y)$ is defined as
    \begin{align*}
        \Hitting(x, y) &= \ep{T_{y} \,|\, X_0 = x}.
    \end{align*}
We recall the following classical fact concerning the relationship between vertex effective resistance and so-called \emph{commute times}, namely, that
    \begin{align}\label{eq:commute-time-identity}
        \Resis(x, y) &= \frac{\Hitting(x, y) + \Hitting(y, x)}{\Vol_w(G)}.
    \end{align}
We refer the reader to, e.g.,~\cite{tetali1991random,vonluxburg2014hitting} for proofs of this fact. It is therefore natural to wonder to what extent~\cref{eq:commute-time-identity} extends to the case of measure effective resistance and, more generally, how to interpret $\Resis(\mu,\nu)$ in terms of vertex hitting times or other related quantities. Na\"{i}vely, any such extension must address how to interpret a hitting time between probability measures in the graph setting. This is facilitated by optimal stopping rules.

\begin{definition}[Stopping rule]
    A \emph{stopping rule} is a map $S$ that associates to each finite path $\tau = (X_0, X_1, \dotsc, X_k)$ on $G$ a number $S(\tau)$ in $[0, 1]$. $S(\tau)$ is understood to be the probability that we continue a random walk given that $\tau$ is the walk thus far observed.
\end{definition}

With slight abuse, we write $\ep{S | X_0\sim \mu}$ to denote the expected length of the corresponding walk conditioned on an initial distribution $\mu\in\Prob(V)$. If $\ep{S | X_0\sim \mu }< \infty$, then the walk stops almost surely in finite time, so we define $X_{S}$ to be the position of the random walk at the stopping time. This leads to a notion of the generalized hitting time between $\mu, \nu$.

\begin{definition}[Access time]\label{def:access-time}
    Let $G=(V,E,w)$ be a connected graph and let $\mu,\nu\in\Prob(V)$. The \emph{access time} $\Hitting(\mu,\nu)$ is given by
        \begin{align}\label{eq:optimal-stopping-rule}
            \Hitting(\mu,\nu) = \inf\left\{ \ep{S \,|\, X_0\sim\mu} : \ep{S |  X_0\sim\mu} < \infty, X_S\sim\nu \right\},
        \end{align}
    where, for any random variable $Y$ on $V$, we say $Y \sim\mu$ if $\pr{Y=x} = \mu(x)$ for each $x\in V$. In other words, $\Hitting(\mu,\nu)$ is the minimum mean length of walks that originate with distribution $\mu$ and terminate according to a stopping rule that achieves distribution $\nu$ at stopping time. If $S^\ast$ achieves the infimum in $\Hitting(\mu,\nu)$, then $S^\ast$ is said to be an \emph{optimal stopping rule}.
\end{definition}

We use the overloaded notation $\Hitting(x, y)$ and $\Hitting(\mu,\nu)$ without risk of confusion since in the case of Dirac measures, the access time recovers the expected hitting time between the corresponding vertices. The set of feasible stopping rules in~\Cref{def:access-time} is nonempty in general provided $G$ is connected. For example, one has the \emph{na\"{i}ve stopping rule} $S_{n}$, obtained as follows: at time $t=0$, perform a sample $x\sim\nu$. Then run $(X_t)_{t\ge 0}$ conditioned on $X_0\sim \mu$, and stop the walk when $X_{S_n} = x$. It follows that $X_{S_{n}} \sim\nu$, and moreover, that
    \begin{align}\label{eq:naive-stopping-rule}
        \mathbb{E}[S_n] = \sum_{x, y\in V} \mu(x)\,\nu(y)\,\Hitting(x,y) < \infty.
    \end{align}
Several examples of optimal stopping rules for general distributions are given in~\cite{lovasz1995mixing}. The authors also show that $\mu,\nu$ always admit an optimal stopping rule achieving the infimum in~\cref{eq:optimal-stopping-rule} as long as $G$ is connected. The access time $\Hitting(\mu,\nu)$ has several useful properties, which we summarize in the proposition below (for proofs, we refer the reader to~\cite{lovasz1995mixing}). 

\begin{proposition}[Properties of access times,~\cite{lovasz1995mixing}]\label{prop:access-time-properties}
    Let $G=(V,E, w)$ be connected and let $\mu,\nu\in\Prob(V)$ be fixed. Then the following facts hold:
        \begin{enumerate}[label=(\roman*)]
            \item If $\nu = \mathbf{1}_y$ for some $y\in V$, then $\Hitting(\mu,\mathbf{1}_y) = \sum_{x\in V}\mu(x) \Hitting(x, y)$,
            \item $\Hitting(\mu,\nu) = \max_{y\in V}\sum_{x\in V} (\mu(x)-\nu(x)) \Hitting(x, y)$,
            \item $\Hitting(\mu,\nu)$ is (separately) convex in its arguments; namely if $\mu,\mu', \nu,\nu'\in\Prob(V)$ are given and $0\le c \le 1$ is fixed, we have
                \begin{align}
                    \Hitting(c \mu + (1-c)\mu', \nu) &\leq c\Hitting(\mu, \nu) + (1-c)\Hitting(\mu', \nu),\text{ and }\\
                    \Hitting(\mu, c \nu + (1-c)\nu' )&\leq c\Hitting(\mu,\nu) + (1-c)\Hitting(\mu, \nu').
                \end{align}
        \end{enumerate}
\end{proposition}

The next result generalizes~\cref{eq:commute-time-identity}.

\begin{theorem}\label{thm:resistance-access-time}
    Let $G=(V,E, w)$ be a connected graph and let $\mu,\nu\in\Prob(V)$ be fixed. Then 
        \begin{align}\label{eq:resistance-access-time}
            \Resis(\mu,\nu) &= -\frac{1}{\mathrm{Vol}_w(G)} \sum_{x, y\in V}(\mu(x) - \nu(x)) (\mu(y) - \nu(y)) \Hitting(x, y).
        \end{align}
\end{theorem}

\begin{proof}
    Let $H\in\mathbb{R}^{V\times V}$ denote the matrix with entries $\Hitting(x, y)$. Then $H - H^\top$ is skew-symmetric, so for any function $f\in\mathbb{R}^V$, it holds that
        \begin{align*}
            f^\top (H - H^\top) f &= 0.
        \end{align*}
    In particular,
        \begin{align*}
            \sum_{x, y\in V}f(x) H(x, y) f(y) &= \frac{1}{2}\sum_{x, y\in V}f(x) (H(x, y) + H(y, x))f(y).
        \end{align*}
    Using~\cref{eq:commute-time-identity}, we have
        \begin{align*}
            \frac{1}{2}\sum_{x, y\in V}f(x) (\Hitting(x, y) + \Hitting(y, x)) f(y)= \frac{\Vol_w(G)}{2}\sum_{x, y\in V}f(x) \Resis(x, y) f(y).
        \end{align*}
    Expanding $\Resis(x, y)$, we have
        \begin{align*}
            \Resis(x, y) &= (\mathbf{1}_x - \mathbf{1}_y)^\top L_{w}^\dagger (\mathbf{1}_x - \mathbf{1}_y) = (L_{w}^\dagger)_{xx} + (L_{w}^\dagger)_{yy} - 2(L_{w}^\dagger)_{xy}.
        \end{align*}
    Thus we obtain
        \begin{align}\label{eq:resis-expanded}
            \frac{\Vol_w(G)}{2}\sum_{x, y\in V}&f(x) \Resis(x, y) f(y) = \frac{\Vol_w(G)}{2}\sum_{x, y\in V}f(x) ((L_{w}^\dagger)_{xx} + (L_{w}^\dagger)_{yy} - 2(L_{w}^\dagger)_{xy}) f(y).
        \end{align}
    With $f=\mu-\nu$, we have $\sum_{x\in V}f(x) = 0$, and thus the diagonal terms in the summand on the right-hand side of~\cref{eq:resis-expanded} vanish, leaving
        \begin{align*}
            \sum_{x, y\in V} (\mu(x) - \nu(x))\Hitting(x, y) (\mu(y) - \nu(y)) &= -\Vol_w(G) \sum_{x, y\in V} (\mu(x) - \nu(x)) (L_{w}^\dagger)_{xy} (\mu(y) - \nu(y))\\
            &= -\Vol_w(G) \Resis(\mu,\nu),
        \end{align*}    
    by definition. The claim is proved.
\end{proof}

We remark that~\Cref{thm:resistance-access-time} is slightly more involved than the comparatively simpler commute time formula for vertex effective resistance. However, it is still possible to bound $\Resis(\mu,\nu)$ in terms of access times. This is summarized in the following corollary.

\begin{corollary}\label{cor:measure-commute-time-inequalities}
    Let $G=(V,E,w)$ be a connected graph and let $\mu,\nu\in\Prob(V)$ be fixed. Then $\Resis(\mu,\nu)$ satisfies the following two inequalities:
        \begin{align}
            \Resis(\mu,\nu) &\leq \frac{2}{\Vol_w(G)} \max\{\Hitting(\mu,\nu), \Hitting(\nu,\mu)\}\\
            \Resis(\mu,\nu) &\leq \frac{1}{\Vol_w(G)} (\Hitting_n(\mu,\nu) + \Hitting_n(\nu, \mu))
        \end{align}
    where $\Hitting_n(\mu,\nu) = \ep{S_n}$ is the expected duration of the na\"{i}ve stopping rule as in~\cref{eq:naive-stopping-rule}.
\end{corollary}

\begin{proof}
    For the first inequality, we use~\Cref{thm:resistance-access-time} and H\"older's inequality:
        \begin{align*}
            \Resis(\mu,\nu) &= -\frac{1}{\Vol_w(G)} \sum_{x, y\in V}(\mu(x) - \nu(x)) (\mu(y) - \nu(y)) \Hitting(x, y)\\
            &= \frac{1}{\Vol_w(G)} \sum_{y\in V} (\nu(y) - \mu(y)) \sum_{x\in V} (\mu(x) - \nu(x)) \Hitting(x, y)\\
            &\leq \frac{1}{\Vol_w(G)} \|\mu-\nu\|_1 \left\| \sum_{x\in V} (\mu(x) - \nu(x)) \Hitting(x, y)\right\|_\infty\\
            &\leq \frac{2}{\Vol_w(G)} \left\| \sum_{x\in V} (\mu(x) - \nu(x)) \Hitting(x, y)\right\|_\infty\\
            &\leq \frac{2}{\Vol_w(G)} \max\{\Hitting(\mu,\nu), \Hitting(\nu,\mu)\},
        \end{align*}
    where the final line follows from~\Cref{prop:access-time-properties}(ii). For the second inequality, we use the vertex commute time formula~\cref{eq:commute-time-identity}. Since $\mu-\nu = \sum_{x,y\in V} \mu(x)\nu(y)(\mathbf{1}_x - \mathbf{1}_y)$ and $z\mapsto z^\top L_w^\dagger z$ is convex, we have
        \begin{align*}
            \Resis(\mu,\nu) &= (\mu-\nu)^\top L_w^\dagger (\mu-\nu)\\
            &\leq \sum_{x,y\in V} \mu(x)\nu(y) (\mathbf{1}_x - \mathbf{1}_y)^\top L_w^\dagger (\mathbf{1}_x - \mathbf{1}_y)\\
            &= \sum_{x,y\in V} \mu(x)\nu(y) \Resis(x, y)\\
            &= \frac{1}{\Vol_w(G)} \sum_{x,y\in V} \mu(x)\nu(y) (\Hitting(x,y) + \Hitting(y,x))\\
            &= \frac{1}{\Vol_w(G)} (\Hitting_n(\mu,\nu) + \Hitting_n(\nu,\mu)),
        \end{align*}
    as claimed.
\end{proof}

Our second combinatorial characterization of resistance distance between measures is related to the spanning trees and spanning $2$-forests of the underlying graph. We are primarily motivated by the following result, originally described in~\cite{bapat1999resistance} (see also~\cite{barrett2020spanning,chung2023forest}).

We denote by $\Tree(G)$ the weighted number of spanning trees of $G$, i.e.,
    \begin{align*}
        \Tree(G) &= \sum_{T\in \mathscr{T}(G)} w_{T}
    \end{align*}
where $\mathscr{T}(G)\subseteq 2^{E}$ is the set of spanning trees of $G$, each viewed as a subset of edges. Similarly, if $x,y\in V$ are fixed, we let $\Forest(x, y)$ denote the weighted number of spanning $2$-forests of $G$ in which $x$ and $y$ belong to different components, defined by
    \begin{align*}
        \Forest(x, y) &= \sum_{F\in \Forest (x, y)} w_{F},
    \end{align*}
where $\Forest (x, y)$ is the set of spanning $2$-forests of $G$ in which $x$ and $y$ belong to different components.

\begin{theorem}[Resistance distance via separating spanning forests,~\cite{bapat1999resistance}]\label{thm:spanning-two-forest}
    Let $G=(V,E,w)$ be a connected graph and let $x,y\in V$ be fixed. Then the effective resistance between $x$ and $y$ admits the representation
        \begin{align*}
            \Resis(x, y) &= \frac{\Forest(x, y)}{\Tree(G)}.
        \end{align*}
\end{theorem}

We now extend~\Cref{thm:spanning-two-forest} to the setting of probability measures on $V$. For $A\subseteq V$ and $\mu\in\Prob(V)$, we write $\mu(A) = \sum_{x\in A} \mu(x)$. For $F\subseteq E$, we write
    \begin{align*}
        \mu(F) := \mu(\{x\in V\,:\, \{x,y\}\in F\text{ for some }y\sim x\}).
    \end{align*}

\begin{theorem}\label{thm:measure-spanning-two-forest}
    Let $G=(V,E,w)$ be a connected graph and let $\mu,\nu\in\Prob(V)$ be fixed. Then 
        \begin{align*}
            \Resis(\mu,\nu) &= \frac{1}{\Tree(G)}\sum_{\substack{F\in\Forest (G) \\ F=T_1\cup T_2}} w_F (\mu(T_1)-\nu(T_1))^2,
        \end{align*}
    where the sum runs over all spanning $2$-forests $F$ of $G$, denoted $\Forest (G)$, with components $T_1$ and $T_2$.
\end{theorem}

\begin{proof}
    We recall from the proof of~\Cref{thm:resistance-access-time}, specifically the argument following~\cref{eq:resis-expanded}, that for any $f\in\mathbb{R}^{V}$ with $\sum_{x\in V}f(x) = 0$, one has
        \begin{align}\label{eq:quad-form-forests}
            f^\top L_w^\dagger f &= -\frac{1}{2}\sum_{x, y\in V} f(x) \Resis(x, y) f(y).
        \end{align}
    By~\Cref{thm:spanning-two-forest}, we have the representation
        \begin{align*}
            \Resis(x, y) &= \frac{1}{\Tree(G)} \sum_{\substack{F\in\Forest (G) \\ F=T_1\cup T_2}} w_{F} \Big(\mathbf{1}\{x\in T_1, y\in T_2\}+\mathbf{1}\{x\in T_2, y\in T_1\}\Big).
        \end{align*}
    where the sum on the right-hand side runs over all spanning $2$-forests of $G$ with components $T_1$ and $T_2$. Applying~\cref{eq:quad-form-forests} with $f = \mu-\nu$, we obtain
        \begin{align*}
            \Resis(\mu,\nu) &= -\frac{1}{2\Tree(G)}\sum_{x, y\in V} \sum_{\substack{F\in\Forest (G) \\ F=T_1\cup T_2}} \Bigl\{  w_F (\mu(x)- \nu(x)) (\mu(y)- \nu(y)) \, \\
            &\qquad\qquad\qquad\qquad\qquad\qquad\qquad\,\times\,\bigl(\mathbf{1}\{x\in T_1, y\in T_2\}+\mathbf{1}\{x\in T_2, y\in T_1\}\bigr) \Bigr\}\\
            &= -\frac{1}{\Tree(G)} \sum_{\substack{F\in\Forest (G) \\ F=T_1\cup T_2}} w_F \sum_{x\in T_1} (\mu(x)- \nu(x))  \sum_{y\in T_2} (\mu(y)- \nu(y)) \\
            &= -\frac{1}{\Tree(G)} \sum_{\substack{F\in\Forest (G) \\ F=T_1\cup T_2}} w_F (\mu(T_1) - \nu(T_1))  (\mu(T_2) - \nu(T_2)) \\
            &= \frac{1}{\Tree(G)} \sum_{\substack{F\in\Forest (G) \\ F=T_1\cup T_2}} w_F (\mu(T_1) - \nu(T_1))^2,
        \end{align*} 
    where the last equality follows from the fact that $\mu(T_1) + \mu(T_2) = \nu(T_1) + \nu(T_2) = 1$, so $\mu(T_2) - \nu(T_2) = -(\mu(T_1) - \nu(T_1))$. The claim is proved.
\end{proof}

We state as a corollary the special case of~\Cref{thm:measure-spanning-two-forest} in the case where the graph $G$ is a tree.

\begin{corollary}\label{cor:resistance-tree-forests}
    Let $G=(V,E,w)$ be a tree, let $E'$ be an orientation of $E$, and let $\mu,\nu\in\Prob(V)$ be fixed. Then
        \begin{align*}
            \Resis(\mu,\nu) &= \sum_{e=(x,y)\in E'} \frac{|\mu(T_e) - \nu(T_e)|^2}{w_e},
        \end{align*}
    where $T_e\subseteq V$ denotes the vertices of the component of $G\setminus e$ that contains the tail of $e$.
\end{corollary}

\begin{remark} 
    \Cref{thm:measure-spanning-two-forest} admits the following probabilistic interpretation. Let $\mathbb{F}$ denote the law on $\Forest (G)$ given by
        \begin{align*}
            \mathbb{F}(F) &\propto w_F,\quad F\in\Forest (G).
        \end{align*}
    For $F\in \Forest (G)$ and $\mu,\nu\in\Prob(V)$ fixed, define the imbalance variable $Z$ according to
        \begin{align*}
            Z(F=T_1\cup T_2) &= \mu(T_1) - \nu(T_1) = \nu(T_2) - \mu(T_2).
        \end{align*}
    Then if $F\sim\mathbb{F}$, $Z$ defines a random variable, and in turn, $\Resis(\mu,\nu)$ is proportional to its second moment:
        \begin{align*}
            \Resis(\mu,\nu) &\propto \mathbb{E}_{F\sim\mathbb{F}}\left({|Z|^2}\right),
        \end{align*}
    where the constant of proportionality depends only on $G$.
\end{remark}

\section{Gradient flows and convexity on the resistance manifold}\label{sec:gradient-flows}

In~\Cref{sec:characterizations}, we introduced the resistance manifold $(\Prob_\star(V),d_\Resis)$ and provided a number of characterizations of $\Resis(\cdot,\cdot)$, equivalently, the corresponding geodesic distance $d_\Resis(\cdot,\cdot)$. The purpose of this section is to analyze gradient flows of various functionals on $\Prob_\star(V)$ with respect to the resistance metric, and in particular to prove~\Cref{thm:gradient-flow} and~\Cref{thm:modulus}.

In this setting, if $F:\Prob_\star(V)\to\mathbb{R}$ is $\mathcal{C}^{1}$, then the Riemannian gradient vector $\nabla F(\rho)\in\mathbb{R}^V_0$ is uniquely determined by
    \begin{align*}
        \partial F(\rho)[f] = f^\top L_w^\dagger \nabla F(\rho),\qquad f\in\mathbb{R}^V_0,
    \end{align*}
where $\partial F(\rho)[f]$ denotes the directional derivative of $F$ at $\rho$ in direction $f$. Given an interval $I\subseteq\mathbb{R}$, a $\mathcal{C}^{1}$ curve $(\rho_t)_{t\in I}\subseteq\Prob_\star(V)$ is called a \emph{gradient flow trajectory of $F$} if it solves
    \begin{align*}
        \frac{d}{dt}\rho_t = -\nabla F(\rho_t),\qquad t\in I.
    \end{align*}

We now define a class of functionals on $\Prob_\star(V)$ given by differentiable $f$-divergences relative to the stationary distribution.

\begin{definition}[$f$-Divergence]\label{defn:f-divergence}
    Let $G=(V,E,w)$ be a connected graph, and let $\pi$ denote its stationary distribution. Let $f:(0,\infty)\to\mathbb{R}$ be continuously differentiable. The (normalized) \emph{$f$-divergence functional} on $\Prob_\star(V)$ is defined by
        \begin{align}\label{eq:f-divergence-functional}
            \Div_f(\rho\,\|\,\pi)
            := \frac{1}{\Vol_w(G)}\sum_{x\in V} \pi(x)\, f\!\left(\frac{\rho(x)}{\pi(x)}\right),
            \qquad \rho\in\Prob_\star(V).
        \end{align}
\end{definition}

The prefactor $1/\Vol_w(G)$ is included for scaling convenience in our analysis of the associated gradient flow. Our next result characterizes the gradient flow of $\Div_f(\cdot\,\|\,\pi)$ on the resistance manifold.

\begin{theorem}\label{thm:gradient-flow-fdiv}
    Let $G=(V,E,w)$ be a connected graph and let $\pi\in\Prob_\star(V)$ denote the stationary distribution of the simple random walk with transition matrix $P=D_w^{-1}A_w$. Let $f:(0,\infty)\to\mathbb{R}$ be $\mathcal{C}^1$, and consider the functional $\Div_f(\cdot\,\|\,\pi)$ defined in~\cref{eq:f-divergence-functional}. Then its gradient is given by
        \begin{align}\label{eq:fdiv-gradient}
            \nabla \Div_f(\rho\,\|\,\pi)
            = \frac{1}{\Vol_w(G)}\, L_w\, f'\!\left(\rho/\pi\right),\quad \rho\in\Prob_\star(V),
        \end{align}
    where $f'(\rho/\pi)\in\mathbb{R}^{V}$ denotes the vector with entries $x\mapsto f'\!\left(\frac{\rho(x)}{\pi(x)}\right)$. Consequently, for any interval $I\subseteq\mathbb{R}$, a smooth curve $(\rho_t)_{t\in I}\subseteq\Prob_\star(V)$ is a gradient flow trajectory of $\Div_f(\cdot\,\|\,\pi)$ if and only if it satisfies the evolution equation
        \begin{align}\label{eq:fdiv-gradient-flow}
            \frac{d}{dt}\rho_t
            = -\frac{1}{\Vol_w(G)}\, L_w\, f'\!\left(\rho_t/\pi\right)
            = \frac{1}{\Vol_w(G)}\, D_w\,(P-I)\, f'\!\left(\rho_t/\pi\right),
        \end{align}
    or, equivalently, for each $x\in V$,
        \begin{align}\label{eq:fdiv-gradient-flow-componentwise}
            \frac{d}{dt}\rho_t(x)
            = \frac{1}{\Vol_w(G)}\sum_{y\sim x} w_{xy}\left(
            f'\!\left(\frac{\rho_t(y)}{\pi(y)}\right) - f'\!\left(\frac{\rho_t(x)}{\pi(x)}\right)\right).
        \end{align}
\end{theorem}

\begin{proof}
    Fix $\rho\in\Prob_\star(V)$ and a tangent direction $h\in\mathbb{R}^V_0$. Differentiating~\cref{eq:f-divergence-functional} at $\rho$ in direction $h$, and using the chain rule, yields
        \begin{align}\label{eq:fdiv-first-variation}
            \partial\Div_f(\rho\,\|\,\pi)[h]
            &= \frac{1}{\Vol_w(G)}\sum_{x\in V}\pi(x)\,f'\!\left(\frac{\rho(x)}{\pi(x)}\right)\frac{h(x)}{\pi(x)}
            = \frac{1}{\Vol_w(G)}\sum_{x\in V} h(x)\, f'\!\left(\frac{\rho(x)}{\pi(x)}\right) \notag\\
            &= \frac{1}{\Vol_w(G)}\, \Bigl(f'(\rho/\pi)\Bigr)^\top h.
        \end{align}
    Let $\Pi_0:\mathbb{R}^V\to\mathbb{R}^V_0$ denote the orthogonal projection onto $\mathbb{R}^{V}_0$. Since $h\in\mathbb{R}^V_0$,~\cref{eq:fdiv-first-variation} can be rewritten as
        \begin{align*}
            \partial\Div_f(\rho\,\|\,\pi)[h]
            = \frac{1}{\Vol_w(G)}\, h^\top \Pi_0 f'(\rho/\pi).
        \end{align*}
    By definition of the resistance gradient, $g:=\nabla \Div_f(\rho\,\|\,\pi)\in\mathbb{R}^V_0$ must satisfy
        \begin{align}\label{eq:fdiv-gradient-identification}
            h^\top L_w^\dagger g = \partial\Div_f(\rho\,\|\,\pi)[h],\qquad \forall\,h\in\mathbb{R}^V_0.
        \end{align}
    Using $L_w^\dagger L_w = L_w L_w^\dagger = \Pi_0$ on $\mathbb{R}^V$, the identities above imply
        \begin{align*}
            g = \frac{1}{\Vol_w(G)}\, L_w\,\Pi_0 f'\!\left(\rho/\pi\right)
            = \frac{1}{\Vol_w(G)}\, L_w\, f'\!\left(\rho/\pi\right),
        \end{align*}
    where the last identity uses that $L_w$ annihilates constants. This proves~\cref{eq:fdiv-gradient}. Substituting this gradient into the defining identity for gradient flow trajectories gives~\cref{eq:fdiv-gradient-flow}, and conversely any smooth curve satisfying~\cref{eq:fdiv-gradient-flow} is a gradient flow trajectory by definition. Expanding $L_w$ as $D_w-A_w$ yields~\cref{eq:fdiv-gradient-flow-componentwise}.
\end{proof}

\Cref{thm:gradient-flow-fdiv} admits two natural corollaries which amount to specializations to the cases of the $\chi^2$ and entropy functionals obtained by different choices of $f$. We state them below without proof.

\begin{corollary}\label{cor:chi2-fdiv}
    Let $f:(0,\infty)\to\mathbb{R}$ be given by $f(r)=\frac{1}{2}(r-1)^2$. Then, for each $\rho\in\Prob_\star(V)$,
        \begin{align*}
            \Div_f(\rho\,\|\,\pi)
            = \frac{1}{2\,\Vol_w(G)}\sum_{x\in V} \frac{(\rho(x)-\pi(x))^2}{\pi(x)}
            =: \chi^2(\rho),
        \end{align*}
    and the gradient flow of $\Div_f(\cdot\,\|\,\pi)$ on $(\Prob_\star(V),d_\Resis)$ satisfies
        \begin{align*}
            \frac{d}{dt}\rho_t
            = -\frac{1}{\Vol_w(G)}\,L_w\left(\frac{\rho_t}{\pi}-\mathbf{1}\right)
            = (P^\top-I)\rho_t.
        \end{align*}
    In particular, this is the continuous time random walk forward equation in column vector form.
\end{corollary}

\begin{corollary}\label{cor:entropy-fdiv}
    Let $f:(0,\infty)\to\mathbb{R}$ be given by $f(r)=r\log r$. Then $\Div_f(\cdot\,\|\,\pi)$ recovers the (normalized) relative entropy functional; specifically, for each $\rho\in\Prob_\star(V)$, we have
        \begin{align*}
            \Div_f(\rho\,\|\,\pi)
            = \frac{1}{\Vol_w(G)}\sum_{x\in V} \rho(x)\log\frac{\rho(x)}{\pi(x)}.
        \end{align*}
    Moreover, the gradient flow of $\Div_f(\cdot\,\|\,\pi)$ on $(\Prob_\star(V),d_\Resis)$ satisfies
        \begin{align*}
            \frac{d}{dt}\rho_t
            = -\frac{1}{\Vol_w(G)}\,L_w\,\log\left(\rho_t/\pi\right),
        \end{align*}
    or, equivalently, for each $x\in V$,
        \begin{align*}
            \frac{d}{dt}\rho_t(x)
            = \frac{1}{\Vol_w(G)}\sum_{y\sim x} w_{xy}\left(
            \log\frac{\rho_t(y)}{\pi(y)}-\log\frac{\rho_t(x)}{\pi(x)}\right).
        \end{align*}
\end{corollary}

We now prove a second result identifying the geodesic strong convexity modulus of $\chi^2$ in the resistance manifold with the spectral gap of the normalized Laplacian. Recall that $\mathcal{L}_w = D_w^{-1/2}L_w D_w^{-1/2}$ is symmetric positive semidefinite and if $G$ is connected, $\mathcal{L}_{w}$ has kernel spanned by $D_w^{1/2}\mathbf{1}$. We write $\lambda_2(\mathcal{L}_w)$ for its smallest strictly positive eigenvalue.

\begin{theorem}\label{thm:modulus-sec4}
    Let $G=(V,E,w)$ be a connected graph and let $\pi$ be the stationary distribution of the simple random walk on $V$. On the resistance manifold $(\Prob_\star(V),d_\Resis)$, the functional $\chi^2$ defined in~\cref{eq:chi2-functional-sec4} has geodesic strong convexity modulus (in the sense of~\cref{eq:geodesic-convexity}) equal to the spectral gap $\lambda_2(\mathcal{L}_w)$.
\end{theorem}

\begin{proof}
    By~\Cref{lem:resistance-geodesics}, every constant-speed geodesic in the resistance manifold is a line segment in $\mathcal{H}$. Fix $\rho_0,\rho_1\in\Prob_\star(V)$, and set $\rho_t=(1-t)\rho_0+t\rho_1$. Let $h:=\rho_1-\rho_0\in\mathbb{R}^V_0$. Then
    using $\pi(x)=\Deg_x/\Vol_w(G)$,
        \begin{align*}
            \chi^2(\rho) = \frac{1}{2}(\rho-\pi)^\top D_w^{-1}(\rho-\pi).
        \end{align*}
    Also, $\rho_t-\pi=(\rho_0-\pi)+t h$, and expanding the quadratic along this segment gives the identity
        \begin{align*}
            \chi^2(\rho_t)
            = (1-t)\chi^2(\rho_0)+t\chi^2(\rho_1) - \frac{t(1-t)}{2}\, h^\top D_w^{-1} h.
        \end{align*}
    By~\Cref{lem:resistance-geodesics}, the squared geodesic distance on the resistance manifold satisfies $d_\Resis(\rho_0,\rho_1)^2 = h^\top L_w^\dagger h$, so $\chi^2$ is geodesic $\lambda$-strongly convex if and only if
        \begin{align*}
            h^\top D_w^{-1} h \ge \lambda\, h^\top L_w^\dagger h,\qquad \forall\,h\in\mathbb{R}^V_0.
        \end{align*}
    Consequently, the largest admissible $\lambda$ is
        \begin{align}\label{eq:lambda-as-rayleigh}
            \lambda_\star := \inf_{h\in\mathbb{R}^V_0\setminus\{0\}} \frac{h^\top D_w^{-1} h}{h^\top L_w^\dagger h}.
        \end{align}

    It remains to identify $\lambda_\star$ with $\lambda_2(\mathcal{L}_w)$. Since $G$ is connected, $\ker(L_w)=\mathrm{span}\{\mathbf{1}\}$ and hence $L_w$ is positive definite on $\mathbb{R}^V_0$. In particular, the map $h\mapsto L_w^{\dagger/2}h$ is a bijection from $\mathbb{R}^V_0$. With the change of variables $z=L_w^{\dagger/2}h$, we have
        \begin{align*}
            h^\top L_w^\dagger h &= \|L_w^{\dagger/2}h\|_2^2 = \|z\|_2^2 = z^\top z,\\
            h^\top D_w^{-1}h &= (L_w^{1/2}z)^\top D_w^{-1}(L_w^{1/2}z) = z^\top\bigl(L_w^{1/2}D_w^{-1}L_w^{1/2}\bigr)z.
        \end{align*}
    Therefore~\cref{eq:lambda-as-rayleigh} becomes
        \begin{align}\label{eq:lambda-rayleigh-LDL}
            \lambda_\star
            = \inf_{z\in\mathbb{R}^V_0\setminus\{0\}}
            \frac{z^\top\bigl(L_w^{1/2}D_w^{-1}L_w^{1/2}\bigr)z}{z^\top z},
        \end{align}
    i.e., $\lambda_\star$ is the smallest strictly positive eigenvalue of the symmetric matrix $L_w^{1/2}D_w^{-1}L_w^{1/2}$. Now define $S:=D_w^{-1/2}L_w^{1/2}$. Then
        \begin{align*}
            S^\top S = L_w^{1/2}D_w^{-1}L_w^{1/2}
            \qquad\text{and}\qquad
            SS^\top = D_w^{-1/2}L_wD_w^{-1/2} = \mathcal{L}_w.
        \end{align*}
    The matrices $S^\top S$ and $SS^\top$ have the same multiset of nonzero eigenvalues. Since $G$ is connected, both have a one-dimensional kernel (spanned by $\mathbf{1}$ for $S^\top S$ and by $D_w^{1/2}\mathbf{1}$ for $\mathcal{L}_w$), and hence their smallest strictly positive eigenvalues agree:
        \begin{align*}
            \lambda_{2}\bigl(S^\top S\bigr) = \lambda_2(SS^\top) = \lambda_2(\mathcal{L}_w).
        \end{align*}
    Combining this with~\cref{eq:lambda-rayleigh-LDL} yields $\lambda_\star=\lambda_2(\mathcal{L}_w)$, as claimed.
\end{proof}

\section*{Acknowledgements}\label{sec:acknowledgements}

The authors would like to acknowledge Yusu Wang for her inspirational discussions early on in this project, Stefan Steinerberger for his help placing the preliminary results within the context of the broader OT space, Zachary Lubberts for his helpful discussions on applications of graph OT to discrete geometry, and Caroline Moosm\"{u}ller for her helpful feedback on linearized Wasserstein distance and Sobolev transport in the continuous setting. SR wishes to acknowledge Evangelos Nikitopoulos, as well as financial support from the Hal\i{}c\i{}o\u{g}lu Data Science Institute through their Graduate Prize Fellowship.  AC was supported by NSF DMS 2012266,  and a gift from Intel, and AC and SR were supported by NSF CISE 2403452. ZW was supported by NSF CCF 2217058.

\end{document}